\documentclass[10pt]{amsart}
\usepackage[francais]{babel}
\usepackage{amsmath,amsfonts}
\usepackage{amssymb}
\usepackage{amscd}
\usepackage{graphics}
\usepackage{graphicx}
\usepackage[all]{xy}
\usepackage{pstcol}

\newtheorem{theoreme}{Th{\'e}or{\`e}me}[section]
\newtheorem{theoremedefinition}[theoreme]{Th{\'e}or{\`e}me et d{\'e}finition}
\newtheorem{proposition}[theoreme]{Proposition}
\newtheorem{lemme}[theoreme]{Lemme}
\newtheorem{corollaire}[theoreme]{Corollaire}
\newtheorem{definition}[theoreme]{Definition}

\newtheorem{remarque}[theoreme]{Remarque}

\def\tr{\mathop{\mathbf{tr}}}
\def\det{\mathop{\mathbf{det}}}

\def\tq{\hskip0.5mm \mid \hskip0.5mm} 

\begin{document}

\title[Primitive du cocycle de Maslov g{\'e}n{\'e}ralis{\'e}]{Primitive du cocycle de Maslov g{\'e}n{\'e}ralis{\'e}}
\author{Jean-Louis Clerc \& Khalid Koufany}
\address{
Institut {\'E}lie Cartan, UMR 7502,
 Universit{\'e} Henri Poincar{\'e} (Nancy 1)
B.P. 239, 54506 Vand{\oe}uvre-l{\`e}s-Nancy cedex, France
}

\email{Jean-Louis.Clerc@iecn.u-nancy.fr, Khalid.Koufany@iecn.u-nancy.fr}
\subjclass{32M15; 53D12; 15A66}
\keywords{Espace hermitien sym{\'e}trique de type tube;
  Fronti{\`e}re de Shilov; Indice d'Arnold-Leray-Maslov; Indice d'inertie;
  Indice de Maslov; Indice de Souriau}

\begin{abstract}
Soit $\mathcal{D}$ un espace hermitien sym{\'e}trique
de type tube, de fronti{\`e}re de Shilov $S$. Nous d{\'e}crivons une
r{\'e}alisation du rev{\^e}tement universel $\widetilde{S}$ de $S$. Nous
construisons  ensuite sur $\widetilde{S}$ une primitive  du cocycle de Maslov
g{\'e}n{\'e}ralis{\'e} tel que d{\'e}fini dans [{\it Transform. Groups} {\bf 6} (2001),
303-320] et [{\it J. Math. Pures Appl.} {\bf 83} (2004), 99-114]. C'est l'analogue de
l'indice de Souriau pour la vari{\'e}t{\'e} lagrangienne. Une variante de cette
construction g{\'e}n{\'e}ralise l'indice d'Arnold-Leray-Maslov. Enfin, nous utilisons cette
primitive pour g{\'e}n{\'e}raliser la notion de nombre de rotation symplectique.
\end{abstract}

\maketitle

\section{Introduction}
Soit $(E,\omega)$ un espace symplectique r{\'e}el, de dimension $2r$. Un {\it
  lagrangien\/} est un sous-espace totalement isotrope (pour la forme $\omega$) de dimension maximale
({\'e}gale {\`a} $r$). L'ensemble des lagrangiens de $E$, not{\'e} $S$, est naturellement muni d'une structure de vari{\'e}t{\'e} compacte et est
appel{\'e} la {\it vari{\'e}t{\'e} lagrangienne\/}. Le groupe symplectique $G=Sp(E)$ op{\`e}re transitivement sur $S$.

A un triplet de lagrangiens $(l_1,l_2,l_3)$, on associe son {\it indice de Maslov\/} $\iota(l_1,l_2,l_3)$, dont une d{\'e}finition alg{\'e}brique  a {\'e}t{\'e}
propos{\'e}e par Kashiwara (voir \cite{Lion-Vergne}). Les
principales propri{\'e}t{\'e}s de l'indice de Maslov sont les suivantes :
\begin{enumerate}
\item[$(i)$] l'indice de Maslov est invariant par $G$, c'est-{\`a}-dire :
$$  \iota(gl_1,gl_2,gl_3) =\iota(l_1,l_2,l_3)$$
pour tout $g\in G$ et $l_1,l_2,l_3\in S$.
\item[$(ii)$] l'indice de Maslov est antisym{\'e}trique pour les permutations des
  trois arguments, c'est-{\`a}-dire :
$$ \iota(l_{\tau(1)},l_{\tau(2)},l_{\tau(3)}) =\epsilon(\tau) \iota(l_1,l_2,l_3)$$
o{\`u} $ \epsilon(\tau)$ d{\'e}signe la signature de la permutation $\tau$.
\item[$(iii)$] l'indice de Maslov satisfait une {\it relation de cocycle\/}
$$ \iota(l_1,l_2,l_3) = \iota(l_1,l_2,l_4) +\iota (l_2,l_3,l_4)+\iota(l_3,l_1,l_4)$$
pour tous $l_1,l_2,l_3,l_4\in S$.\\
\end{enumerate}

La vari{\'e}t{\'e} lagrangienne est un cas particulier d'une famille d'espaces, li{\'e}e
aux {\it alg{\`e}bres de Jordan euclidiennes\/}. On renvoie {\`a} \cite{Faraut-Koranyi} pour les
d{\'e}finitions, notations et r{\'e}sultats concernant ces alg{\`e}bres.\\

Selon un premier point de vue, on voit $S$ comme la {\it fronti{\`e}re de Shilov
d'un domaine de type tube\/}. Plus pr{\'e}cis{\'e}ment, soit
$J$ une alg{\`e}bre de Jordan euclidienne,
$\Omega$ son c{\^o}ne (ouvert, propre, convexe, homog{\`e}ne, auto-adjoint) des
carr{\'e}s inversibles. Notons $\mathbb{J}$ l'alg{\`e}bre de Jordan complexifi{\'e}e de
$J$. Elle est naturellement munie d'une {\it norme spectrale\/}. La
boule unit{\'e} ouverte correspondante
$$\mathcal{D} = \{z\in \mathbb{J} \; , \; |z|<1\}$$
munie de la m{\'e}trique de Bergmann, est un domaine hermitien sym{\'e}trique
sous l'action du groupe $ G$ des diff{\'e}omorphismes holomorphes de
$\mathcal{D}$. Par une transform{\'e}e
de Cayley, le domaine $\mathcal{D}$ est holomorphiquement {\'e}quivalent  au tube
$$ T_\Omega = J+i\Omega$$
d'o{\`u} le nom de domaine de type tube.
On d{\'e}signe par
$S$ la fronti{\`e}re de Shilov de $\mathcal{D}$. Lorsque
$J$ est l'alg{\`e}bre $Sym(r,\mathbb{R})$ des matrices $r\times r$ r{\'e}elles sym{\'e}triques
pour le produit $ x\cdot y = \frac{1}{2}(xy+yx)$, le domaine $\mathcal{D}$ correspondant
est le {\it disque de Siegel}, le groupe $G$ est essentiellement le groupe
symplectique et la fronti{\`e}re de Shilov s'identifie {\`a} la vari{\'e}t{\'e}
lagrangienne {\'e}voqu{\'e}e ci-dessus.

Il existe un deuxi{\`e}me point de vue, qui ne fait pas appel {\`a} la th{\'e}orie des
domaines hermitiens sym{\'e}triques. La vari{\'e}t{\'e} $S$ peut {\^e}tre obtenue comme {\it
compl{\'e}t{\'e}e conforme\/} de l'alg{\`e}bre de Jordan
$J$, le groupe $G$ apparaissant comme le groupe des {\it transformations conformes}
de $J$ suivant la terminologie de \cite{Bertram}. Une variante de ce point de vue
consiste {\`a} munir la vari{\'e}t{\'e}
$S$ d'une structure causale, et
$G$ est alors le groupe des transformations causales de $S$.\\

Dans \cite{Clerc2}, achevant le programme commenc{\'e} dans \cite{Clerc-Orsted-1} et
\cite{Clerc-Orsted-2}, l'un des deux pr{\'e}sents auteurs a propos{\'e} une
g{\'e}n{\'e}ralisation de l'indice de Maslov pour l'action de $G$ sur $S\times S\times S$. L'approche,
nouvelle m{\^e}me pour le cas de la vari{\'e}t{\'e} lagrangienne, utilise principalement la
r{\'e}alisation de
$S$ comme fronti{\`e}re de Shilov du domaine $\mathcal{D}$. L'objectif du pr{\'e}sent
article est de faire le lien avec les constructions classiques dues {\`a} Maslov,
Arnold, Souriau et Leray (voir \cite{Maslov}, \cite{Arnold},
\cite{Arnold2}, \cite{Souriau}, \cite{Leray}, \cite{deGosson2}). Plus pr{\'e}cis{\'e}ment, l'indice de Maslov
satisfait une relation de cocycle. Lorsqu'on passe au rev{\^e}tement universel
$\widetilde{S}$ de $S$, on cherche une {\it primitive\/} de ce cocycle, c'est-{\`a}-dire
une fonction $m:\widetilde{S}\times \widetilde{S}\longrightarrow \mathbb{Z}$, antisym{\'e}trique et
telle que
$$ \iota(\sigma_1,\sigma_2,\sigma_3) = m(\widetilde{\sigma}_1, \widetilde{\sigma}_2)+m(\widetilde{\sigma}_2,\widetilde{\sigma}_3)+m(\widetilde{\sigma}_3,\widetilde{\sigma}_1)$$
pour tout triplet $(\widetilde{\sigma}_1, \widetilde{\sigma}_2,\widetilde{\sigma}_3)$ de points de
$\widetilde{S}$, de projections respectives $\sigma_1,\sigma_2,\sigma_3$ sur $S$.
Dans une variante, on introduit l'{\it indice d'inertie\/} d'un triplet de
points de $S$. Celui-ci poss{\`e}de encore la propri{\'e}t{\'e} d'invariance par le groupe
$G$ et une propri{\'e}t{\'e} de cocycle (mais n'a pas la propri{\'e}t{\'e}
d'antisym{\'e}trie), et on se propose {\'e}galement d'en chercher une primitive.

Les deux constructions, quoique essentiellement inspir{\'e}es par les travaux
ant{\'e}rieurs, en diff{\`e}rent sur quelques points. L'utilisation de la structure
causale sur $S$ et la description d'un rev{\^e}tement $\Gamma$ du groupe
$G$ qui op{\`e}re sur $\widetilde{S}$ sont nouvelles, m{\^e}me dans le cas classique, et
permettent notamment de donner des d{\'e}finitions dont l'invariance par $\Gamma$
est manifeste. \\

Le plan de l'article est le suivant.

Le paragraphe 2 pr{\'e}sente les domaines de
type tube et leurs fronti{\`e}res de Shilov, du point de vue des alg{\`e}bres de Jordan
euclidiennes. La structure causale sur $S$ y est {\'e}galement introduite.
 Au
paragraphe 3, on donne une r{\'e}alisation du rev{\^e}tement universel $\widetilde{S}$ de $S$,
puis on construit un rev{\^e}tement $\Gamma$ du groupe
$G$ qui op{\`e}re sur $\widetilde{S}$.
 Au paragraphe 4, on {\'e}tudie l'action de $G$ sur
$S\times S$, et on pr{\'e}cise un invariant caract{\'e}ristique pour les
couples de points de $S$, appel{\'e} {\it indice de transversalit{\'e}\/}.

Apr{\`e}s
ce travail pr{\'e}paratoire, on commence au paragraphe 5 la construction de
l'{\it indice de Souriau\/} dans le cas transverse, {\'e}tendu au cas non
transverse au paragraphe 6, en suivant une id{\'e}e de de Gosson (voir \cite{deGosson}).

Au paragraphe 7,
on introduit le {\it cycle de Maslov\/} et on {\'e}tablit ses
propri{\'e}t{\'e}s g{\'e}om{\'e}triques. Elles permettent au paragraphe 8 de d{\'e}finir le
{\it nombre d'Arnold\/} (directement  sans hypoth{\`e}se de transversalit{\'e}). On
construit ensuite l'{\it indice d'Arnold-Leray-Maslov\/}, qui fournit une
primitive de l'{\it indice d'inertie\/} (d{\'e}fini par modification de l'indice de
Maslov).

Au paragraphe 9, on montre, sans
s'appesantir sur les d{\'e}tails, que le concept de l'{\it indice de Maslov de deux
chemins\/} tel que d{\'e}velopp{\'e} dans \cite{Cappell-Lee-Miller} peut aussi se g{\'e}n{\'e}raliser dans
notre contexte. Enfin, au paragraphe 10, on construit pour le groupe $G$
l'analogue du {\it nombre de rotation symplectique\/} introduit dans \cite{Barge-Ghys}.

\section{La fronti{\`e}re de Shilov d'un espace hermitien sym{\'e}trique de type tube}

Soit $J$ une alg{\`e}bre de Jordan euclidienne simple d'{\'e}l{\'e}ment unit{\'e} $e$, de
dimension $n$ et de rang $r$. Son {\it invariant de Peirce} $d$ est tel que
\begin{eqnarray*}
  n=r+\frac{r(r-1)}{2}d.
\end{eqnarray*}
Soit $\det$ le d{\'e}terminant de $J$. C'est une fonction polynomiale
homog{\`e}ne de degr{\'e} $r$. La trace de $J$, not{\'e}e $\tr$ est une forme
lin{\'e}aire sur $J$. On peut supposer que le
produit salaire de $J$ est d{\'e}fini par
\begin{equation}\label{ps}
\langle x|y \rangle=\tr(xy).
\end{equation}
Par complexification on obtient une alg{\`e}bre de Jordan complexe
$\mathbb{J}$. La conjugaison par rapport {\`a} la forme r{\'e}elle $J$ est
not{\'e}e $\eta$, mais on utilisera aussi la notation
\begin{equation}\label{conjugaison}
\eta(z)=\overline{z}.
\end{equation}
 Les extensions holomorphes de $\det$ et $\tr$ {\`a} $\mathbb{J}$ seront
 encore not{\'e}es $\det$ et $\tr$ respectivement.\\
Le produit scalaire (\ref{ps}) s'{\'e}tend en un produit hermitien sur $\mathbb{J}$ donn{\'e} par
\begin{eqnarray}
  \label{fh}
  \langle z|w \rangle=\mathbf{tr}(z\overline{w}).
\end{eqnarray}
Soit 
\begin{equation*}
S=\{\sigma\in\mathbb{J} \tq \overline{\sigma}=\sigma^{-1}\}.
\end{equation*}
C'est une sous-vari{\'e}t{\'e} compacte connexe de $\mathbb{J}$. Nous allons
r{\'e}sumer certaines de ses caract{\'e}risations dans la proposition suivante
(voir \cite[Proposition X.2.3]{Faraut-Koranyi}) : 
\begin{proposition} \label{proprietes_S}
Pour $\sigma \in\mathbb{J}$, les assertions suivantes sont
  {\'e}quivalentes : 
  \begin{enumerate}
  \item[$(i)$] $\sigma\in S$,
  \item[$(ii)$] $\sigma=\exp(iu)$ o{\`u} $u\in J$,
  \item[$(iii)$] il existe un rep{\`e}re de Jordan $(c_j)_{1\leq j\leq r}$ de $J$ et des
    nombres complexes $\zeta_1,\ldots,\zeta_r$ de module 1, tels que
    $\sigma=\sum_{j=1}^r\zeta_jc_j$.\\
  \end{enumerate}
\end{proposition}

On note $\mathbb{L}=Str(\mathbb{J})$ le groupe de structure de l'alg{\`e}bre
de Jordan complexe $\mathbb{J}$. Consid{\'e}rons
\begin{equation*}
L(S)=\{g\in GL(\mathbb{J}) \tq g(S)=S\}.
\end{equation*}
 Alors $L(S)=\mathbb{L}\cap
U(\mathbb{J})$, o{\`u} $U(\mathbb{J})$ d{\'e}signe le groupe unitaire de
$\mathbb{J}$ muni de la forme hermitienne (\ref{fh}). Ce groupe op{\`e}re
transitivement sur $S$. Le stabilisateur de $e$ dans $L(S)$ co{\"\i}ncide avec le
groupe $Aut(J)$ des automorphismes de l'alg{\`e}bre de Jordan $J$ (on {\'e}tend
l'action des {\'e}l{\'e}ments de $Aut(J)$ de mani{\`e}re $\mathbb{C}-$lin{\'e}aire
{\`a} $\mathbb{J}$). L'involution $\eta$ introduite en (\ref{conjugaison})
conserve $S$, et admet $e$ comme point fixe isol{\'e}. L'involution correspondante de
$L(S)$ donn{\'e}e par $ g\mapsto \eta \circ g \circ \eta$ a pour
points fixes exactement $L(S)_e = Aut(J)$. En d'autres termes, muni de
la m{\'e}trique induite par le produit hermitien (\ref{fh}),
$S$ est un espace riemannien sym{\'e}trique compact, isomorphe {\`a}
$L(S)/Aut(J)$.\\
Comme $L(S)$ n'est pas n{\'e}cessairement connexe, on pr{\'e}f{\`e}re introduire
sa composante connexe neutre $U$, qui est encore transitive sur $S$, et on
d{\'e}signe par $U_e$ le stabilisateur dans $U$ de l'{\'e}l{\'e}ment $e$. On a
$Aut(J)^\circ\subset U_e \subset Aut(J)$ et 
\begin{equation*}
S\simeq U/U_e.
\end{equation*}
 
 Soit $(c_j)_{1\leq j\leq r}$ un rep{\`e}re de Jordan. Pour tout {\'e}l{\'e}ment $z\in
\mathbb{J}$, il existe $u\in U$ et il existe des scalaires
$0\leq\lambda_1\leq \ldots\leq\lambda_r$ uniques tels que
\begin{equation*}
z=u\Bigl(\sum_{j=1}^r\lambda_jc_j\Bigr).
\end{equation*}
La {\it norme spectrale} de $z$ est d{\'e}finie par
\begin{equation*}
|z|=\sup\lambda_j
\end{equation*} 
et la fonction $z\mapsto |z|$ est une norme sur $\mathbb{J}$ invariante par
$U$. 
Soit $\mathcal{D}$ la boule unit{\'e} ouverte de $\mathbb{J}$ pour la norme spectrale
\begin{equation*}
\mathcal{D}=\{z\in \mathbb{J}\;,\; |z|<1\}.
\end{equation*}

\begin{theoreme}[{\rm voir \cite[Theorem X.4.6]{Faraut-Koranyi}}]
$\mathcal{D}$ est un domaine born{\'e} hermitien de type tube dont la fronti{\`e}re de
Shilov est $S$. 
\end{theoreme}
Le domaine $\mathcal{D}$ peut {\^e}tre r{\'e}alis{\'e} comme un domaine tube au dessus d'un
c{\^o}ne ouvert de $J$. Plus pr{\'e}cis{\'e}ment, l'ensemble des carr{\'e}s de $J$ est un c{\^o}ne convexe ferm{\'e} et son int{\'e}rieur $\Omega$ est un c{\^o}ne ouvert,
convexe, propre, et auto-dual. C'est le {\it c{\^o}ne sym{\'e}trique} associ{\'e}
{\`a} $J$.\\
Le groupe $L(\Omega)$ des transformations lin{\'e}aires de $J$ qui pr{\'e}servent
$\Omega$ est un groupe de Lie r{\'e}ductif qui agit  transitivement sur
$\Omega$. Soit $L$ sa composante connexe neutre. Le
stabilisateur  $K$ de l'identit{\'e} $e$ dans $L$ est un sous-groupe compact
maximal de $L$. C'est la composante connexe neutre du groupe $Aut(J)$ et on a $K=L\cap O(J)$, o{\`u}
$O(J)$ est le groupe orthogonal de $J$.

Soit $T_\Omega$ le {\it domaine tube} dans $\mathbb{J}$ de base $\Omega$
\begin{equation*}
  T_\Omega=J+i\Omega=\{z=x+iy \tq x\in J,\; y\in \Omega\}.
\end{equation*}
La {\it transformation de Cayley} $c$ et son inverse $p$ sont donn{\'e}es par 
\begin{eqnarray*}
  \label{eq:cayley}
  p(z)&=(z-ie)(z+ie)^{-1}&=e-2i(z+ie)^{-1},\\
  c(w)&=i(e+w)(e-w)^{-1}&=-ie+2i(e-w)^{-1},
\end{eqnarray*}
et on a (voir  \cite[Proposition X.2.3, Theorem X.4.3]{Faraut-Koranyi}) :
\begin{theoreme} L'application $p$ induit un isomorphisme biholomorphe de
  $T_\Omega$ sur $\mathcal{D}$, et
\begin{equation*}
p(J)=\{\sigma\in S \tq \det(e-\sigma)\not=0\}.\\
\end{equation*}
\end{theoreme}
Les deux domaines $\mathcal{D}$ et $T_\Omega$ sont holomorphiquement
{\'e}quivalents. $T_\Omega$ est la r{\'e}alisation non born{\'e}e de $\mathcal{D}$ et
on peut voir $J$ comme sa fronti{\`e}re de Shilov; son image par $p$ est un ouvert dense dans $S$.\\

Soit $G$ la composante connexe neutre du groupe des
automorphismes holomorphes de $\mathcal{D}$. C'est un groupe de Lie
semi-simple et le stabilisateur de $0$ dans $G$ est un sous-groupe compact
maximal de $G$ et co{\"\i}ncide avec le groupe $U$. 

Soit $G(T_\Omega)$ la composante connexe neutre du groupe des
automorphismes holomorphes de $T_\Omega$. Les {\'e}l{\'e}ments du groupe $L$  op{\`e}rent sur $T_\Omega$ par extension
$\mathbb{C}-$lin{\'e}aire et on
peut consid{\'e}rer $L$ comme  un sous-groupe de $G(T_\Omega)$. En fait
$G(T_\Omega)$ est engendr{\'e} par $L$, le sous-groupe ab{\'e}lien $N^+$ des
translations 
\begin{equation}\label{inversion}
t_u:z\mapsto z+u,\;\; u\in J
\end{equation}
et par l'inversion 
\begin{equation*}
s : z\mapsto -z^{-1}.
\end{equation*}
Le sous-groupe ab{\'e}lien $N^-=s\circ N^+\circ s$ est l'ensemble des
transformations
\begin{equation}\label{transmutation}
\widetilde{t}_u=s\circ t_u\circ s : z\mapsto (z^{-1}+u)^{-1},\;\; u\in J.
\end{equation}
On a $$G(T_\Omega)=c\circ G\circ c^{-1}.$$

Le produit semi-direct $P=L\ltimes N^-$ est un
sous-groupe parabolique maximal de $G(T_\Omega)$. L'espace
homog{\`e}ne $G(T_\Omega)/P$ est donc une vari{\'e}t{\'e} (r{\'e}elle) compacte contenant $J$
comme un sous-ensemble ouvert dense, le plongement de $J$ dans
$G(T_\Omega)/P$ {\'e}tant donn{\'e} par 
\begin{equation}\label{plongement_de_J}
J\to G(T_\Omega)/P :  u\mapsto t_uP.
\end{equation}
$G(T_\Omega)/P$ est la {\it compactification  conforme} de $J$ (voir \cite{Bertram}), elle est isomorphe de la fronti{\`e}re
de Shilov $S$ de $\mathcal{D}$.\\

Donnons enfin un troisi{\`e}me point de vue sur la g{\'e}om{\'e}trie de $S$, {\`a}
savoir l'existence d'une structure causale $G-$invariante sur $S$.

Un {\it c{\^o}ne causal \/} d'un espace vectoriel de dimension finie est
un c{\^o}ne non vide $C$, convexe, ouvert et d'adh{\'e}rence pointue (c'est-{\`a}-dire $\overline{C}\cap -\overline{C}=\{0\}$).
Une {\it structure causale\/} sur une vari{\'e}t{\'e} $\mathcal{M}$ est la donn{\'e}e
d'un champs de c{\^o}nes $(C_p)_{p\in\mathcal{M}}$ sur $\mathcal{M}$
tel que, pour tout $p\in \mathcal{M}$, $C_p$ est un c{\^o}ne causal dans
l'espace tangent $T_p(\mathcal{M})$, qui d{\'e}pend diff{\'e}rentiablement
de $p$.
 Une courbe $\gamma : [a,b]\to \mathcal{M}$ de classe $\mathcal{C}^1$
est {\it une courbe causale} (resp. {\it anti-causale}) si, pour tout $t$,
$\dot{\gamma}(t)\in C_{\gamma(t)}$ (resp. $\dot{\gamma}(t)\in
-C_{\gamma(t)}$).\\
Lorsque  $\mathcal{M}$ est un espace homog{\`e}ne sous l'action d'un groupe de
Lie $G$, la
structure causale est dite {\it $G-$invariante} si
$C_{g(x)}=Dg(x)(C_x)$, pour tout $g\in G$ et tout $x\in\mathcal{M}$.  Soit $x_\circ$ un point base de $\mathcal{M}$. Une structure causale
$G-$invariante sur $\mathcal{M}$ est compl{\`e}tement d{\'e}termin{\'e}e par la donn{\'e}e d'un
c{\^o}ne causal $C$, dans $T_{x_\circ}(\mathcal{M})$, invariant par l'action du
stabilisateur de $x_\circ$ dans  $G$.

L'alg{\`e}bre de Jordan $J$ admet une structure causale naturelle,
model{\'e}e par le c{\^o}ne sym{\'e}trique $\Omega$. Il s'agit de la donn{\'e}e
en chaque point de $J$ du c{\^o}ne $\Omega$, vu comme un c{\^o}ne
dans l'espace tangent {\`a} $J$ en ce point. Cette structure causale est
invariante par $G(T_\Omega)$.

La transformation de Cayley permet de transporter {\`a} $S$ cette structure
causale. Consid{\'e}rons $-e\in S$ comme point
base de $S$, alors l'espace tangent {\`a} $S$ en $-e$ s'identifie {\`a}
$iJ$. Par ailleurs, $p(0)=-e$, et la transformation de Cayley $p=c^{-1}$ est bien d{\'e}finie
dans un voisinage de $0$. Sa diff{\'e}rentielle  en $0$ est
$$Dp(0)=-2i\text{Id}_{\mathbb{J}}.$$ 
On d{\'e}finit alors un c{\^o}ne causal en
$-e$ en posant 
$$C_{-e}=Dp(0)(\Omega)=-i\Omega.$$
Comme ce c{\^o}ne est invariant par le stabilisateur du point $-e$ dans $G$,
on peut d{\'e}finir une structure causale sur $S$ invariante par $G$ de la
mani{\`e}re suivante : 
 si $\sigma$
est un point quelconque de $S$, alors il existe $g\in G$ tel que
$g(-e)=\sigma$. Le c{\^o}ne
$$C_\sigma=Dg(-e)(C_{-e}),$$
est un c{\^o}ne causal contenu dans l'espace tangent {\`a} $S$ en $\sigma$, ne
d{\'e}pendant pas du choix de $g$. Nous
obtenons ainsi un champs de c{\^o}nes causaux $(C_\sigma)_{\sigma\in S}$ invariants par
$G$. 
Enfin, il est clair que le c{\^o}ne causal $C_\sigma$ d{\'e}pend
diff{\'e}rentiablement de $\sigma$.

Pour finir ce paragraphe, nous donnons la liste des alg{\`e}bres de Jordan
euclidiennes, les domaines hermitiens born{\'e}s correspondants et leurs
fronti{\`e}res de Shilov.

\begin{table}[h]
\begin{tabular}[t]{cccc}
\hline
$J$  & $\mathbb{J}$ & $\mathcal{D}$ & $S$ \\
\hline
$Sym(m,\mathbb{R})$ & $Sym(m,\mathbb{C})$ & $Sp(2m,\mathbb{R})/U(m)$ & ${U(m)/O(m)}$ \\
${Herm(m,\mathbb{C})}$ & ${Mat(m,\mathbb{C})}$ & ${SU(m,m)/S(U(m)\times U(m))}$& ${U(m)}$\\
${Herm(m,\mathbb{H})}$ & ${Skew(2m,\mathbb{C})}$ & ${SO^*(4m)/U(2m)}$ & ${U(2m)/SU(m,\mathbb{H})}$\\
${\mathbb{R}\times\mathbb{R}^{q-1}}$ & ${\mathbb{C}\times\mathbb{C}^{q-1}}$ &
${SO_0(2,q)/SO(2)\times SO(q)}$ & ${(U(1)\times S^{q-1})/\mathbb{Z}_2}$\\
${Herm(3,\mathbb{O})}$ & ${Mat(3,\mathbb{O})}$ & ${E_{7(-25)}/U(1)E_6}$ &
${U(1)E_6/F_4}$\\
\hline
\end{tabular}
\caption{Domaines hermitiens born{\'e}s et leurs fronti{\`e}res de Shilov}
\end{table}

\begin{table}[h]
\begin{center}
\begin{tabular}[t]{cccc}
\hline
${J}$  & ${n}$ & ${r}$ & ${d}$\\
\hline
${Sym(m,\mathbb{R})}$                 & ${\frac{1}{2}m(m+1)}$ & ${m}$ & ${1}$\\
${Herm(m,\mathbb{C})}$               & ${m^2}$               & ${m}$ & ${2}$  \\
${Herm(m,\mathbb{H})}$               & ${m(2m-1)}$           & ${m}$ & ${4}$\\
${\mathbb{R}\times\mathbb{R}^{q-1}}$ & ${q}$                 & ${2}$   & ${q-2}$\\
${Herm(3,\mathbb{O})}$               & ${27}$                  & ${3}$   &${8}$\\
\hline
\end{tabular}
\caption{La dimension, le rang et l'invariant de Peirce de $J$}
\end{center}
\end{table}

\section{Le rev{\^e}tement universel de la fronti{\`e}re de Shilov}
La vari{\'e}t{\'e} $S$ vue comme espace riemannien sym{\'e}trique $\simeq U/U_e$  n'est
pas semi-simple mais seulement r{\'e}ductif. En
effet, $S$ est stable par la multiplication par un nombre complexe
$e^{i\theta}$ de module $1$, et ceci montre qu'il y a un sous-groupe
isomorphe au groupe des nombres complexes de module $1$ dans le centre de $U$.\\
Consid{\'e}rons maintenant
\begin{equation*}
S_\mathtt{1}=\{\sigma\in S \tq \det(\sigma)=1\}.
\end{equation*}

\begin{proposition}
$S_\mathtt{1}$ est une sous-vari{\'e}t{\'e} connexe de $S$. C'est un espace riemannien sym{\'e}trique semi-simple de type compact.
\end{proposition}
\proof
Pr{\`e}s du point $e\in S$, on peut utiliser la carte locale de $S$ donn{\'e}e
par
$$J\ni x \longmapsto \exp ix\in S.$$
La condition $\sigma\in S_\mathtt{1}$ se traduit par $\tr(x)=0$, ce qui
montre bien que $S_\mathtt{1}$ est une sous-vari{\'e}t{\'e} pr{\`e}s de $e$. On utilise
alors l'action du groupe $U$ pour en d{\'e}duire que $S_\mathtt{1}$ est
(globalement) une sous-vari{\'e}t{\'e} de $S$. Pour la connexit{\'e}, on va montrer
que $S_\mathtt{1} = \exp iJ_\mathtt{0} $, o{\`u} 
\begin{equation*}
  J_\mathtt{0}=\{v\in J\tq \tr(v)=0\}.
\end{equation*}
Comme $\det(\exp iv) = e^{i\tr(v)}=1$ pour tout $v\in J_\mathtt{0}$, seule
l'inclusion
$ S_\mathtt{1}\subset\exp iJ_\mathtt{0}$ est {\`a} d{\'e}montrer. Soit donc $\sigma\in S_\mathtt{1}$.
D'apr{\`e}s la proposition \ref{proprietes_S} on a $S=\exp iJ$, il existe donc $v\in J$ tel que $\exp iv=\sigma$. Comme
$$1=\det(\sigma) = \det(\exp iv) = e^{i\tr(v)}, $$ on a $\tr(v) = 2k\pi$ avec
$k\in \mathbb{Z}$. Soit $(c_j)_{1\leq j\leq r}$ un rep{\`e}re de Jordan
adapt{\'e} {\`a} $v$, c'est-{\`a}-dire tel que $ v=\sum_{j=1}^r a_j c_j$, o{\`u}
$a_j\in\mathbb{R}$ pour tout $j$. Soit
$$w=v-2k\pi c_1=(a_1-2k\pi)c_1+\sum_{j=2}^r a_j c_j.$$ Alors 
$$\tr(w) = \tr(v) - 2k\pi\tr(c_1)=0,$$ car $\tr(c)=1$ pour tout idempotent primitif $c$, et $$\exp
w= e^{i(a_1-2k\pi)}c_1 + \sum_{j=2}^r e^{ia_j}c_j= \sum_{j=1}^r
e^{ia_j}c_j= \exp v = \sigma\ .$$ Donc $\sigma\in \exp i J_\mathtt{0}$ et par
suite $S_\mathtt{1} = \exp iJ_\mathtt{0}$, ce qui montre que $S_\mathtt{1}$ est connexe.\\
Le polyn{\^o}me $\det$ est
relativement invariant par l'action du groupe de structure $\mathbb{L}$ de $\mathbb{J}$, et il y a donc un
caract{\`e}re $\chi$ de $\mathbb{L}$ {\`a} valeurs dans $\mathbb{C}^*$ tel que
\begin{equation*}
  \det(gz)=\chi(g)\det(z)
\end{equation*}
pour tout $g\in \mathbb{L}$ et tout $z\in \mathbb{J}$. Soit 
\begin{equation*}
U_\mathtt{1} = \{u\in U\mid \chi(u) = 1\} .
\end{equation*} 
Le sous-groupe $U_\mathtt{1}$ op{\`e}re sur $S_\mathtt{1}$,  et cette action est 
transitive. En effet, si $\sigma\in S_\mathtt{1}$, il existe $u\in U$ tel que $\sigma
=u(e)$. Alors $$1=\det(\sigma) = \chi(u) \det(e)=\chi(u),$$ et donc $u\in
U_\mathtt{1}$. \\
La composante connexe neutre $U_\mathtt{1}^\circ$ de $U_\mathtt{1}$ op{\`e}re donc aussi
transitivement sur $S_1$. De plus, si $u\in U$ est un {\'e}l{\'e}ment qui fixe $e$,
alors $\chi(u) = 1$ et donc $U_e \subset U_\mathtt{1}$. En particulier, $(U_\mathtt{1})_e $ est le
sous-groupe des points fixes de l'involution  $g\mapsto \eta\circ g\circ
\eta$ (o{\`u} $\eta$ est donn{\'e}e par (\ref{conjugaison})), ce qui r{\'e}alise $S_\mathtt{1}$ comme espace sym{\'e}trique 
\begin{equation*}
S_\mathtt{1}=U_\mathtt{1}^\circ/(U_\mathtt{1}^\circ\cap Aut(J)).
\end{equation*}
 Reste {\`a} voir que $U_\mathtt{1}^\circ$ est semi-simple.
L'alg{\`e}bre de Lie de $U$ est $\mathfrak{u} = \mathfrak{k}\oplus i\mathfrak{p}$, o{\`u}
$\mathfrak{k} = Der(J)$ est l'alg{\`e}bre des d{\'e}rivations de l'alg{\`e}bre de Jordan
$J$, et 
$$\mathfrak{p} = \{ L(v)\tq v\in J\},$$ o{\`u} $L(v)$ d{\'e}signe l'op{\'e}rateur de
multiplication par $v$ dans $J$. Donc clairement
$$\mathfrak{u}_1 = \mathrm{Lie}(U_\mathtt{1})= \mathfrak{k} \oplus i\mathfrak{p}_0$$
o{\`u} 
$$\mathfrak{p}_0 = \{L(v)\tq v\in J, \tr(v) = 0\}=\{L(v) \tq v\in J_{\mathtt{0}}\}.$$ Le centre de $\mathfrak{u}$ est
$\mathbb{R} L(e)=\mathbb{R}\mathrm{Id}$, de sorte que le centre de $\mathfrak{u}_1$ est
r{\'e}duit {\`a} $\{0\}$, et donc $\mathfrak{u}_1$ est bien semi-simple.
\hfill\hfill $\square$ \\

Pour d{\'e}cider si un espace sym{\'e}trique semi-simple de type compact $U/K$ est
simplement connexe (ou plus g{\'e}n{\'e}ralement pour calculer son groupe
fondamental), rappelons qu'il faut comparer deux r{\'e}seaux dans un
sous-espace de Cartan  : l'un , not{\'e} $\Lambda_0$ est obtenu
alg{\'e}briquement {\`a} partir de donn{\'e}es radicielles, l'autre $\Lambda$ est
obtenue {\`a} partir de la d{\'e}termination des {\'e}l{\'e}ments $H$ d'un sous-espace
de Cartan tels que $\exp H\in K$. 

Fixons un rep{\`e}re de Jordan $(c_j)_{1\leq j\leq r}$, et soit
\begin{equation*}
\mathfrak{a} = \{L(a)\tq a=\sum_{j=1}^r a_j c_j,\; a_j\in \mathbb{R}\}.
\end{equation*}
C'est un sous-espace de Cartan de $\mathfrak{p}$, et par suite 
\begin{equation*}
 \mathfrak{a}_0 = \{L(a)\tq a=\sum_{j=1}^r a_j c_j,\; a_j\in
\mathbb{R},\;\tr(a) = 0\}.
\end{equation*}
est un sous-espace de Cartan de $\mathfrak{p}_0$. On munit $\mathfrak{p}$ (et
donc $\mathfrak{a}$) du produit scalaire d{\'e}fini par 
$$(L(x)|L(y))=\langle x|y\rangle,$$
 pour $x, y\in J$. Il est proportionnel au produit scalaire induit par la forme de Killing de $\mathfrak{u}$. \\
 Le syst{\`e}me des racines
restreintes de $(\mathfrak{u},\mathfrak{a})$  est 
$$\Sigma(\mathfrak{u},\mathfrak{a})=\{{\frac{1}{2} }(a_j-a_k),1\leq j\neq k\leq
r\},$$ (voir \cite[page 212]{Faraut-Koranyi}) et c'est aussi le syst{\`e}me des racines
restreintes de $(\mathfrak{u}_0, \mathfrak{a}_0)$. \\
Consid{\'e}rons les deux r{\'e}seaux de $\mathfrak{a}_0$ :
\begin{equation*}
\Lambda_0=\text{r{\'e}seau engendr{\'e} par }\{2\pi\frac{A_{j,k}}{(
  A_{j,k}|A_{j,k} )} \tq 1\leq j\not= k\leq r\},
\end{equation*}
o{\`u} $A_{j,k}$ est le covecteur de $\frac{1}{2}(a_j-a_k)$, et
\begin{equation*}
\Lambda=\{H\in \mathfrak{a}_0 \tq \exp(iH)e=e\}.
\end{equation*}

\begin{lemme}\label{lemme_Lambda_0} Le r{\'e}seau $\Lambda_0$ est engendr{\'e} par les {\'e}l{\'e}ments
\begin{equation*}
2\pi L(c_j-c_k),\quad 1\leq j\not= k\leq r.
\end{equation*}
\end{lemme}
\proof
En effet, l'{\'e}l{\'e}ment $A_{j,k}=\frac{1}{2}L(c_j-c_k)\in \mathfrak{a}_0$ est le
covecteur \footnote{On devrait utiliser le produit scalaire
induit par la forme de Killing, mais tout produit scalaire sur $\mathfrak{a}_0$ qui lui est proportionnel convient dans la d{\'e}termination du r{\'e}seau
$\Lambda_0$.}  de la racine
${\frac{1}{2}}(a_j-a_k)$ en ce sens que pour $a=\sum_{j=1}^ra_jc_j$,
$$ \langle \frac{1}{2}(c_j-c_k)|a \rangle = \frac{1}{2}(a_j-a_k).$$
Notons par ailleurs que 
$$( A_{j,k}|A_{j,k} ) = \langle
\frac{1}{2}(c_j-c_k)|\frac{1}{2}(c_j-c_k) \rangle =
\frac{1}{2}\ .$$ D'o{\`u}
$$2\pi \frac{A_{j,k}}{( A_{j,k}|A_{j,k} )} =
2\pi L(c_j-c_k).$$
\hfill\hfill$\square$

\begin{lemme}\label{lemme_Lambda} On a
\begin{equation*} 
\Lambda=\{ L(a)\tq a=\sum_{j=1}^r 2\pi m_j c_j,\, m_j\in \mathbb{Z},\,\sum_{j=1}^r m_j=0\}.
\end{equation*}
\end{lemme}
\proof
Soit $a=\sum_{j=1}^ra_jc_j$. On a 
$$ \exp iL(a) = P(\exp(\frac{ia}{2})),$$ o{\`u}
$P$ d{\'e}signe l'application quadratique de l'alg{\`e}bre de Jordan $J$, et $$P(\exp  (\frac{ia}{2}))e =\exp(ia) .$$
Donc $$L(a)\in \Lambda\Longleftrightarrow  \exp(ia)=e\ .$$
Comme $ \exp(ia) = \sum_{j=1}^r e^{ia_j}c_j $, cette condition {\'e}quivaut
encore {\`a} $ a_j\in 2\pi\mathbb{Z}$ pour tout entier $j$ tel que $1\leq j\leq
r$. D'o{\`u} le lemme. \hfill\hfill$\square$

\begin{theoreme}\label{S_1_SC}
L'espace sym{\'e}trique $S_\mathtt{1}$ est simplement connexe.
\end{theoreme}
\proof Les lemmes \ref{lemme_Lambda_0} et \ref{lemme_Lambda}
montrent que les r{\'e}seaux $\Lambda_0$ et $\Lambda$ co{\"\i}ncident. Donc d'apr{\`e}s
\cite[Theorem 3.6]{Loos} (voir aussi \cite[Ch. VII. Theorem 8.4 et Theorem 9.1]{Helgason2}), l'espace
$S_\mathtt{1}$ est bien simplement connexe. \hfill\hfill$\square$\\

Suivant une m{\'e}thode classique, nous r{\'e}alisons maintenant le
rev{\^e}tement universel de $S$. Soit
\begin{equation*}
\widetilde{S} = \{(\sigma,\theta)\in S\times \mathbb{R}\tq \det(\sigma)= e^{ir\theta}\},
\end{equation*}
muni de la topologie induite par la topologie produit sur $S\times \mathbb{R}$. La projection sur le premier facteur r{\'e}alise $\widetilde{S}$ comme un
rev{\^e}tement de $S$ avec fibre isomorphe {\`a} $\mathbb{Z}$.

\begin{theoreme} 
$\widetilde{S}$ est le rev{\^e}tement universel de $S$.
\end{theoreme}
\proof 
Montrons en effet que l'application 
\begin{equation*}
S_\mathtt{1}\times \mathbb{R}\longrightarrow \widetilde{S}\quad (\sigma, \theta)\longmapsto (e^{i\theta}\sigma, \theta)
\end{equation*}
est un hom{\'e}omorphisme. D'abord l'application est bien {\`a} valeurs dans
$\widetilde{S}$, car $\det(e^{i\theta}\sigma)= e^{ir\theta}\det(\sigma)=e^{ir\theta}$. Ensuite, l'application est clairement bijective,
l'application r{\'e}ciproque {\'e}tant
$(\sigma,\theta)\longmapsto(e^{-i\theta}\sigma, \theta)$. De plus la
continuit{\'e} de l'application et de sa r{\'e}ciproque est imm{\'e}diate. Comme $S_\mathtt{1}$
est simplement connexe d'apr{\`e}s le th{\'e}or{\`e}me \ref{S_1_SC}, l'espace
$\widetilde{S}$ l'est aussi. \hfill\hfill$\square$\\

Nous allons maintenant construire un rev{\^e}tement $\Gamma$ de $G$,
puis une op{\'e}ration de $\Gamma$ sur
$\widetilde{S}$ relevant l'action de $G$ sur $S$.\\

 Soit $g\in G$ et $z\in\mathcal{D}$. On note
$J(g,z) = Dg(z)$ la diff{\'e}rentielle  de la transformation $
z\longmapsto g(z)$. On sait que c'est un {\'e}l{\'e}ment de $\mathbb{L}$, et
on pose 
$$j(g,z) =
\chi(J(g,z))\in \mathbb{C}^*.$$
 Pour $g\in G$ fix{\'e}, on consid{\`e}re les
d{\'e}terminations de l'argument de $j(g,\cdot)$. Plus pr{\'e}cis{\'e}ment une
fonction $ \varphi(g,\cdot) : \mathcal{D}\longrightarrow \mathbb{R}$ est
appel{\'e}e une {\it d{\'e}termination de l'argument $j(g,\cdot)$\/} de si elle est
continue et si 
$$\forall z\in \mathcal{D},\;\; e^{i\varphi(g,z)}=\frac{j(g,z)}{|j(g,z)|}.$$ L'existence de telles
d{\'e}terminations est cons{\'e}quence du fait que
$\mathcal{D}$ est simplement connexe. Deux telles d{\'e}terminations diff{\`e}rent d'un
multiple entier de $2\pi$.\\
On pose alors
$$ \Gamma = \bigl\{\big(g,\varphi_g\big)\;|\; g\in G\bigr\} $$
o{\`u} $g\in G$, et $\varphi_g = \varphi(g,\cdot)$ est une d{\'e}termination de
l'argument de $j(g,\cdot)$.  La loi de groupe est donn{\'e}e par
$$ (g,\varphi(g,\cdot))(h,\psi(h,\cdot))=\big(gh, \varphi(g,h(\cdot))+\psi(h,\cdot)\big).$$
Remarquons en effet que la formule
$$ J(gh,z) = J(g,h(z))J(h,z)$$
implique que 
$$
e^{i\left[\varphi(g,h(z))+\psi(h,z)\right]}=e^{i\varphi(g,h(z))}e^{i\psi(h,z)}
=\frac{j(g,h(z))}{|j(g,h(z))|}\frac{j(h,z)}{|j(h,z)|}
= \frac{j(gh,z)}{|j(gh,z)|}.
$$
Par suite, $\varphi(g,h(\cdot))+\psi(h,\cdot)$ est bien une d{\'e}termination de
l'argument de $ j(gh,\cdot)$.\\

Pour la topologie, on remarque qu'une d{\'e}termination $\varphi(g,\cdot)$ de
l'argument de
$j(g,\cdot)$ est d{\'e}termin{\'e}e par  sa valeur $\varphi(g,0)$ {\`a} l'origine. Donc
$\Gamma$ s'identifie ensemblistement au sous-ensemble ferm{\'e} de
$G\times\mathbb{R}$ donn{\'e} par
$$ \big\{(g,\theta)\in G\times \mathbb{R}\;\mid\; e^{i\theta}= j(g,0)\big\},$$ et on munit $\Gamma$ de la topologie induite, faisant ainsi
de $\Gamma$ un groupe topologique (la r{\'e}f{\'e}rence \cite{Godement} contient
les d{\'e}tails n{\'e}cessaires dans le cas o{\`u} $ G=PSL(2,\mathbb{R})$). 
La projection $\pi:\Gamma
\longrightarrow G$ est continue, et fait de $\Gamma$ un rev{\^e}tement de
$G$\, \footnote{On observera que le groupe $\Gamma$ n'est pas n{\'e}cessairement
  connexe, comme on peut le voir dans le cas o{\`u} $G=Sp(E)$.}.

L'application $\iota : \mathbb{Z}\longrightarrow \Gamma$ donn{\'e}e par
$\iota: n\longmapsto (\text{id},2n\pi)$
est un homomorphisme de groupes, et la suite exacte 
\begin{equation*}
 \begin{CD}
  0 @> >> \mathbb{Z} @>\iota>> \Gamma @>\pi >> G @> >> 1
 \end{CD} 
\end{equation*}
est une extension centrale de $G$. 

\begin{lemme}
Soit $g\in G$, et $\sigma\in S$. Alors
\begin{equation}\label{equation_det_j}
\det(g(\sigma)) = \frac{j(g,\sigma)}{|j(g,\sigma)|} \det(\sigma).
\end{equation}
\end{lemme}
\proof Soit $G_\sigma$ le stabilisateur de $\sigma$ dans
$G$. Comme $U$ est transitif sur $S$, on a $G=UG_\sigma$. L'{\'e}l{\'e}ment $g$
peut donc s'{\'e}crire $g=up$, avec $u\in U$ et $p\in G_\sigma$. Le membre de
gauche  de la formule (\ref{equation_det_j}) vaut
$$\det(g(\sigma))=\det(u(\sigma))=\chi(u)\ \det(\sigma) .$$
D'autre part, $j(g,\sigma)=j(u,\sigma)j(p,\sigma)$. Maintenant
$J(p,\sigma)$ conserve globalement le plan tangent {\`a} $S$ en $\sigma$, qui
est une forme r{\'e}elle de $\mathbb{J}$, et conserve globalement le c{\^o}ne
$C_\sigma$, et donc $ j(p,\sigma)$ est un nombre r{\'e}el strictement
positif. D'autre par, l'action de $u$ {\'e}tant lin{\'e}aire, on a pour tout $z\in
\mathbb{J}$, $ J(u,z)=u$, et $j(u,z) = \chi(u)$. Par suite,
$$\frac{j(g,\sigma)}{|j(g,\sigma)|}=\frac{j(u,\sigma)}{|j(u,\sigma)|}\frac{j(p,\sigma)}{|j(p,\sigma)|}=
\chi(u) .$$
\hfill\hfill$\square$\\

On peut maintenant d{\'e}finir une action de $\Gamma$ sur $\widetilde{S}$.

\begin{proposition}
Pour $ (g,\varphi(g,\cdot))\in \Gamma$
et
$(\sigma,\theta)\in \widetilde{S}$, posons
$$ \big(g,\varphi(g,\cdot)\big)
(\sigma,\theta)=\big(g(\sigma),\theta+\frac{1}{r}\varphi(g,\sigma)\big) .$$
On d{\'e}finit ainsi une action continue de $\Gamma$ sur $\widetilde{S}$.
\end{proposition}
\proof On v{\'e}rifie en effet que
$$e^{i(r\theta + \varphi(g,\sigma))}=e^{i\varphi(g,\sigma)}e^{ir\theta}=\frac{j(g,\sigma)}{|j(g,\sigma)|}\det(\sigma) =\det(g(\sigma))
$$
d'apr{\`e}s le lemme pr{\'e}c{\'e}dent, et donc
$(g(\sigma),\theta+\frac{1}{r}\varphi(g,\sigma))$ est bien {\'e}l{\'e}ment de
$\widetilde{S}$. 
Ensuite
pour $(g,\varphi(g,\cdot)) $ et $ (h,\psi(h,\cdot)$ {\'e}l{\'e}ments de $\Gamma$, on a
$$
\begin{array}{ll}
(g,\varphi(g,\cdot)\Big((h,\psi(h,\cdot))(\sigma,\theta)\Big)&=(g,\varphi(g,\cdot))\big(h(\sigma),\theta+\frac{1}{r}\psi(g,\sigma)\big)\\
&=\big(g(h(\sigma)), \theta +\frac{1}{r}\psi(g,\sigma)+\frac{1}{r}\varphi(g,h(\sigma)\big)\\
&=\big(gh,\psi(g,\cdot)+\varphi(g,h\cdot)\big)(\sigma,\theta)\\
&= \big((g,\varphi(g,\cdot)(h,\psi(h,\cdot)\big)(\sigma,\theta)
\end{array}
$$
ce qui montre que l'on a bien {\`a} faire {\`a} une action (clairement continue) de
$\Gamma$ sur
$\widetilde{S}$. 

\hfill\hfill$\square$

Notons que l'{\'e}l{\'e}ment $T=\iota(1)$ op{\`e}re sur $\widetilde{S}$ par 
\begin{equation}\label{action_par_T}
T\cdot(\sigma,\theta) = (\sigma,\theta+\frac{2\pi}{r}).
\end{equation}
\medskip

\section{L'indice de transversalit{\'e}}
Nous allons introduire dans ce paragraphe un invariant qui
caract{\'e}rise les orbites de l'action de $G$ sur $S\times S$. Transport{\'e} par la
transformation de Cayley, cet invariant  est li{\'e} au rang des
{\'e}l{\'e}ments de l'alg{\`e}bre de Jordan $J$. Cette notion g{\'e}n{\'e}ralise la
{\it distance arithm{\'e}tique\/} {\'e}tudi{\'e}e par Hua \cite{Hua} dans le cadre des
alg{\`e}bres de Jordan $Sym(m,\mathbb{R})$ et $Herm(m,\mathbb{C})$.\\

Rappelons d'abord la caract{\'e}risation des orbites de $L$ sur $J$.
Soit $e=\sum_{j=1}^rc_j$ une d{\'e}composition de Peirce de l'identit{\'e}.
Pour tout couple d'entiers
$0\leq p, q\leq r$, tels que $0\leq p+q\leq r$, on pose
$$e_{p,q}= \sum_{j=1}^p c_j-\sum_{j=p+1}^{p+q}c_j .$$
Notons $\mathcal{O}_{p,q}=L e_{p,q}$ l'orbite du point $e_{p,q}$ sous
l'action du  groupe $L$.
\begin{theoreme}\label{theo_inertie_Sylvester}
L'espace $J$ est r{\'e}union disjointe des orbites
$\mathcal{O}_{p,q}$,  pour $0\leq p,q\leq r$, $p+q\leq r$.
\end{theoreme}
\proof (esquisse; pour une autre d{\'e}monstration, voir
  \cite{Kaneyuki})
Soit $x\in J$. Soient $\lambda_j,1\leq j\leq
r$ ses valeurs propres, figurant autant de fois que leur
multiplicit{\'e}. Soit $p$ (resp. $q$) le nombre de valeurs propres
strictement positives (resp. strictement n{\'e}gatives). D'apr{\`e}s
le th{\'e}or{\`e}me spectral, $x$ est conjugu{\'e} par $K$ {\`a} l'{\'e}l{\'e}ment
$\sum_{j=1}^{p+q} \lambda_jc_j$, o{\`u} on peut supposer de
plus que $\lambda_j>0$ pour $1\leq j\leq p$ et $\lambda_j<0$ pour
$p+1\leq j\leq p+q$. Gr{\^a}ce {\`a} l'action du sous-groupe
$A=\exp L(\mathfrak{a})$, on voit  que l'orbite de $x$ contient
l'{\'e}l{\'e}ment $ e_{p,q}$. Donc $J$ est bien r{\'e}union des orbites $\mathcal
{O}_{p,q} $.

Pour voir qu'elles sont disjointes, on observe que la
signature de la forme quadratique $Q_x: y\longmapsto \langle P(x)y,y\rangle$ ne d{\'e}pend que de l'orbite de $x$. Cela r{\'e}sulte de ce
que $ P(gx)=gP(x)g^t$ pour tout $g\in L$. Il suffit alors de
v{\'e}rifier que les signatures des formes quadratiques $Q_{p,q}$
donn{\'e}es par
$ Q_{p,q}(y) = \langle P(e_{p,q})y,y\rangle$
sont toutes distinctes.

En particulier, il y a $r+1$ orbites ouvertes, correspondant aux
cas $p+q=r$, dont la r{\'e}union est l'ouvert $J^\times$ des {\'e}l{\'e}ments
inversibles de $J$.
\hfill\hfill$\square$\\

Le {\it rang} d'un {\'e}l{\'e}ment $x$ de $J$, not{\'e} $\text{rang}(x)$, est par
d{\'e}finition le nombre de ses valeurs propres non nulles compt{\'e}es avec leur multiplicit{\'e} (voir
\cite[Theorem III.1.1]{Faraut-Koranyi}). C'est un invariant sous l'action
du groupe $L(\Omega)$. On observe que $x$ et $y$ ont m{\^e}me rang si et
seulement si $P(x)$ et $P(y)$ ont m{\^e}me rang.\\

Soit $c$ un idempotent de $J$.  La d{\'e}composition de Peirce de $J$
relativement {\`a} $c$ est (voir \cite[chapter IV]{Faraut-Koranyi})
\begin{equation}\label{decomp_Peirce}
J=J(c,1)\oplus J(c,\frac{1}{2})\oplus J(c,0).
\end{equation}
Si $c$ est de rang $k$, alors 
\begin{eqnarray}\label{dim_J_0}
\dim J(c,1)&=&k+\frac{k(k-1)}{2}d,\nonumber \\
\dim J(c,0)&=&(r-k)+\frac{(r-k)(r-k-1)}{2}d,\\
\dim J(c,\frac{1}{2})&=&k(r-k)d. \nonumber
\end{eqnarray}
Tout {\'e}l{\'e}ment $x\in J$ se d{\'e}compose relativement {\`a} (\ref{decomp_Peirce})
sous la forme
\begin{equation*}
x=x_1+x_{\frac{1}{2}}+x_0 \;\;\mathrm{avec}\;\; x_0\in J(c,0),\;x_1\in J(c,1),\;x_{\frac{1}{2}}\in J(c,\frac{1}{2}).
\end{equation*}
Les
sous-espaces
$J(c,1)$ et $J(c,0)$ sont des sous-alg{\`e}bres de Jordan de
$J$. On d{\'e}signe
par  $\Omega_1$ (resp. $\Omega_0$) le
c{\^o}ne sym{\'e}trique associ{\'e} {\`a} $J(c,1)$ (resp. $J(c,0)$) et on note $\det_c$
(resp. $\det_{e-c}$) le d{\'e}terminant de $J(c,1)$ (resp. $J(c,0)=J(e-c,1)$).\\

Voici quelques lemmes qui nous seront utiles dans la suite.

\begin{lemme}\label{lemme_unicite_c}
Soit $x\in J$ un {\'e}l{\'e}ment de rang $k$. Il existe un unique idempotent $c$
de rang $k$ tel que $cx=x$. De plus $c\in\mathbb{R}[x]$.
\end{lemme}
\proof D'apr{\`e}s le th{\'e}or{\`e}me spectral, (voir \cite[Theorem
III.1.1]{Faraut-Koranyi}), il existe un syst{\`e}me complet d'idempotents
 orthogonaux \footnote{non n{\'e}cessairement primitifs}, $(c_j)_{1\leq j\leq m}$, un  $m-$uplet de nombres r{\'e}els
$(\lambda_j)_{1\leq j\leq m}$ deux {\`a} deux distincts tels que
$x=\sum_{j=1}^m\lambda_jc_j$. Les $(c_j)_{1\leq j\leq m}$ appartiennent {\`a}
$\mathbb{R}[x]$, et les $(c_j,\lambda_j)_{1\leq j\leq m}$ sont uniques {\`a}
l'ordre pr{\`e}s. Le rang de $x$ est
la somme des rangs des idempotents $c_j$ pour lesquels
$\lambda_j\not=0$. L'{\'e}l{\'e}ment
$$c=\sum_{1\leq j\leq m, \;\lambda_j\not=0}c_j$$
est un idempotent de m{\^e}me rang que $x$, et satisfait $cx=x$.\\
R{\'e}ciproquement, soit $d$ un idempotent de rang $k=\mathrm{rang}(x)$, tel
que $dx=x$. Appliquant le th{\'e}or{\`e}me spectral dans l'alg{\`e}bre
$J(d,1)$, on voit qu'il existe une famille unique d'idempotents $(d_j)_{1\leq
  j\leq p}$ deux {\`a} deux orthogonaux, de somme $\sum_{j=1}^pd_j=d$ et des
scalaires uniques $(\mu_j)_{1\leq j\leq p}$ deux {\`a} deux distincts tels
que $x=\sum_{j=1}^d\mu_jd_j$. De plus aucun des $\mu_j$ n'est nul, puisque
$\mathrm{rang}(x)=\mathrm{rang}(d)$ par hypoth{\`e}se. Soit
$d_{p+1}=e-d$. Alors $\sum_{j=1}^{p+1}d_j=e$ et on peut {\'e}crire
$$x=\sum_{j=1}^p\mu_jd_j+0d_{p+1}.$$
D'apr{\`e}s l'assertion d'unicit{\'e} dans le th{\'e}or{\`e}me spectral, on voit
que $m=p+1$, et qu'{\`a} l'ordre pr{\`e}s, les $d_j$ (resp. $\mu_j$)
co{\"\i}ncident avec les $c_j$ (resp. $\lambda_j$). D'o{\`u} facilement $c=d$.
\hfill\hfill$\square$

\begin{lemme}\label{lemme_Omega_1_Omega_0}
 Soit $x=x_1+x_{\frac{1}{2}}+x_0$ la d{\'e}composition de Peirce d'un
 {\'e}l{\'e}ment $x\in J$  relativement {\`a} un idempotent $c\in J$. Si $x\in \Omega$, alors $x_1\in\Omega_1$ et $x_0\in\Omega_0$.
\end{lemme}
\proof Rappelons d'abord une caract{\'e}risation du c{\^o}ne
$\Omega$ (voir \cite[page 4]{Faraut-Koranyi}). Pour tout $x\in J$, on a l'{\'e}quivalence suivante
\begin{equation*}
x\in \Omega \iff \forall y\in\overline{\Omega}\setminus\{0\},\; \langle x|y\rangle>0.
\end{equation*}
On a une caract{\'e}risation analogue pour $\Omega_1\subset J(c,1)$, que nous allons
utiliser pour d{\'e}montrer le lemme. Soit donc
$y_1\in\overline{\Omega}_1\setminus\{0\}$. D'apr{\`e}s le th{\'e}or{\`e}me
spectral relativement {\`a} $J(c,1)$ (voir \cite[Theorem III.1.2]{Faraut-Koranyi}), il
existe une d{\'e}composition de Peirce de l'unit{\'e}
$c=\sum_{j=1}^kc_j$ telle que l'on puisse {\'e}crire
\begin{equation*}
y_1=\sum_{j=1}^k\lambda_jc_j,
\end{equation*}
les $\lambda_j$, $1\leq j\leq k$ {\'e}tant positifs ou nuls et non tous nuls. Il
en r{\'e}sulte en particulier que $y_1\in\overline{\Omega}\setminus\{0\}$. En utilisant
l'orthogonalit{\'e} de la d{\'e}composition (\ref{decomp_Peirce}) on voit que
\begin{equation*}
\langle x_1|y_1\rangle=\langle x_1+x_\frac{1}{2}+x_0|y_1\rangle=\langle x|y_1\rangle>0.
\end{equation*}
Ceci {\'e}tant vrai pour tout $y_1\in\overline{\Omega}_1\setminus\{0\}$, on a
bien $x\in\Omega_1$. Pour montrer que $x_0\in\Omega_0$, il suffit de
consid{\'e}rer  la d{\'e}composition de Peirce par rapport {\`a} l'idempotent $e-c$.

\hfill\hfill$\square$\\

\begin{lemme}\label{lemme_det}
Soit $c$ un idempotent de $J$. Pour tout $x_1\in J(c,1)$ et tout $x_0\in J(c,0)$ on a
\begin{equation*}
\det(x_1+x_0)=\det{}\!_c(x_1)\det{}\!_{e-c}(x_0).
\end{equation*}
\end{lemme}
\proof Soit $k$ le rang  de $c$. Alors $J(c,1)$
et $J(c,0)=J(e-c,1)$ sont deux alg{\`e}bres de Jordan de rang $k$ et $r-k$
respectivement. Il existe alors $\{c_1,\ldots,c_k\}$ et $\{c_{k+1},\ldots,
c_r\}$ deux rep{\`e}res de Jordan de $J(c,1)$
et $J(c,0)$ respectivement tels que $x_1=\sum_{j=1}^k\lambda_jc_j$ et
$x_0=\sum_{j=k+1}^r\lambda_jc_j$ soient les d{\'e}compositions spectrales de
$x_1\in J(c,1)$ et $x_0\in J(c,0)$. Or $\{c_1\ldots,
c_k,c_{k+1},\ldots, c_r\}$ est un rep{\`e}re de Jordan de $J$, donc
$\sum_{j=1}^r\lambda_jc_j$ est la d{\'e}composition spectrale de
$x_1+x_0$ dans $J$. On en d{\'e}duit que 
$$\det(x_1+x_0)=(\prod_{j=1}^k\lambda_j)(\prod_{j=k+1}^r\lambda_j)=\det{}\!_c(x_1)\det{}\!_{e-c}(x_0).$$
\hfill\hfill$\square$\\

Pour tout $z, w\in J$, on introduit l'op{\'e}rateur 
$$z\Box w=L(zw)+[L(z),L(w)].$$
{\`A} chaque $z$ de $J(c,\frac{1}{2})$ est associ{\'e} la {\it
  transformation de Frobenius}
\begin{equation*}\label{Frobenius}
\tau_c(z)= \exp (2z\Box c)
\end{equation*}
 qui est un {\'e}l{\'e}ment du groupe $L$.\\
Si $x=x_1+x_{\frac{1}{2}}+x_0$ est la d{\'e}composition de Peirce de $x\in J$
par rapport {\`a} $c$, alors (voir \cite[Lemma VI.3.1]{Faraut-Koranyi}) $\tau_c(z)x=y_1+y_{\frac{1}{2}}+y_0$, o{\`u}
\begin{eqnarray}\label{decomp_frobenius1}
y_1&=&x_1 \nonumber \\
y_{\frac{1}{2}}&=&2L(z)x_1+x_{\frac{1}{2}}\\
y_0&=&2L(e-c)\left[L(z)^2x_1+L(z)x_{\frac{1}{2}}\right]+x_0 \nonumber
\end{eqnarray}
On en d{\'e}duit que $\mathrm{Det}(\tau_c(z))=1$. Soit $\tau_c(z)^*$
l'adjoint de la transformation de Frobenius $\tau_c(z)$ par rapport au
produit scalaire (\ref{ps}). On a
\begin{equation*}
\tau_c(z)^*=\exp(2c\Box z)=\exp(2z\Box(e-c))=\tau_{e-c}(z).
\end{equation*}
Si $u=u_1+u_{\frac{1}{2}}+u_0$ est la d{\'e}composition de Peirce d'un {\'e}l{\'e}ment $u\in J$ par rapport {\`a} $c$, alors la d{\'e}composition $\tau_c(z)^*u=v_1+v_{\frac{1}{2}}+v_0$ est donn{\'e}e par
\begin{eqnarray}\label{decomp_frobenius2}
v_1&=&2L(c)\left[L(z)^2u_0+L(z)u_{\frac{1}{2}}\right]+u_1 \nonumber \\
v_{\frac{1}{2}}&=&2L(z)u_0+u_{\frac{1}{2}}\\
v_0&=&u_0\nonumber 
\end{eqnarray}
En particulier, la restriction de $\tau_c(z)^*$ {\`a} $J(c,1)$ co{\"\i}ncide avec
l'identit{\'e}.\\

Montrons maintenant le lemme suivant  qui g{\'e}n{\'e}ralise la proposition VI.3.2
dans \cite[page 107]{Faraut-Koranyi} {\`a} un idempotent $c$ quelconque.

\begin{lemme}\label{lemme_Forbenuis}
Soit $c$ un idempotent quelconque de $J$, et soit $J=J(c,1)\oplus J(c,\frac{1}{2})\oplus J(c,0)$ la
d{\'e}composition de $J$ associ{\'e}e. Soit $\omega\in\Omega$ et
soit $\omega=\omega_1+\omega_{\frac{1}{2}}+\omega_0$ sa d{\'e}composition
correspondante. Alors il
  existe un unique $z\in J(c,\frac{1}{2})$  tel que
\begin{equation*}
\omega=\tau_c(z)(\omega_1+\xi_0)
\end{equation*}
o{\`u} $\xi_0$ est un {\'e}l{\'e}ment du c{\^o}ne sym{\'e}trique $\Omega_0$ associ{\'e}
{\`a} l'alg{\`e}bre de Jordan $J(c,0)$.
\end{lemme}
\proof L'alg{\`e}bre de
Jordan $J(c,1)$ agit par multiplication sur $J(c,\frac{1}{2})$, et l'application $\Phi $ de $J(c,1)$ dans $\mathrm{End} (J(c,\frac{1}{2}))$ d{\'e}finie par
\begin{equation}\label{representation_J_c_1}
\forall x\in J(c,1),\; \forall \xi\in J(c,\frac{1}{2}),\;\; \Phi(x)\xi=2x\xi
\end{equation}
 est une repr{\'e}sentation
de l'alg{\`e}bre de Jordan
$J(c,1)$ sur $J(c,\frac{1}{2})$.\\ 
Si $z\in J(c,\frac{1}{2})$ alors d'apr{\`e}s les relations
(\ref{decomp_frobenius1}) 
\begin{eqnarray*}
\tau_c(-z)\omega&=&\omega_1-2L(z)\omega_1+\omega_{\frac{1}{2}}+\xi_0\\
         &=& \omega_1-\Phi(\omega_1)z+\omega_{\frac{1}{2}}+\xi_0,
\end{eqnarray*}
o{\`u} $\xi_0\in J(c,0)$. Du fait que $\omega\in \Omega$, d'apr{\`e}s le lemme \ref{lemme_Omega_1_Omega_0}, la
composante $\omega_1$ est inversible dans $J(c,1)$ et par suite l'endomorphisme $\Phi(\omega_1)$ est
inversible dans $\mathrm{End}(J(c,\frac{1}{2}))$. La proposition s'en d{\'e}duit en choisissant
$z=\Phi(\omega_1)^{-1}\omega_{\frac{1}{2}}$.

\hfill\hfill$\square$

\begin{lemme}\label{lemme_estim_det}
Soit $x\in J$ un {\'e}l{\'e}ment de rang $k$, $c$ l'unique idempotent de rang $k$
tel que $cx=x$, et soit $J=J(c,1)\oplus J(c,\frac{1}{2})\oplus J(c,0)$ la
d{\'e}composition de $J$ associ{\'e}e {\`a} l'idempotent $c$. Soit $\omega\in\Omega$ et
soit $\omega=\omega_1+\omega_{\frac{1}{2}}+\omega_0$ sa d{\'e}composition
correspondante. Alors
\begin{equation*}
\det(x+t\omega)=t^{r-k}\det\!{}_c(x)\det\!{}_{e-c}(\omega_0)+O(t^{r-k+1}).
\end{equation*}
\end{lemme}
\proof
D'apr{\`e}s le lemme \ref{lemme_Forbenuis} il existe $z\in J(c,\frac{1}{2})$
tel que $\omega=\tau_c(z)^*(\omega_0+\xi_1)$, avec $\xi_1\in\Omega_1$. Comme
on a de plus $\tau_c(z)^*x=x$, il vient
$$x+t\omega=\tau_c(z)^*\bigl((x+t\xi_1)+t\omega_0\bigr)$$
et donc d'apr{\`e}s le lemme \ref{lemme_det}
$$\det(x+t\omega)=\det\!{}_c(x+t\xi_1)t^{r-k}\det\!{}_{e-c}(\omega_0)$$
d'o{\`u} le r{\'e}sultat. \hfill\hfill$\square$\\

\begin{corollaire}\label{corollaire_D_det}
Soit $x\in J$ un {\'e}l{\'e}ment de rang $k$ et soit $ \Phi_x$ la fonction d{\'e}finie
sur $J$ par $ \Phi_x(z) = \det (x+z)$. Alors,
\begin{itemize}
\item[$(i)$] pour tout $j$, $1\leq j\leq r-k-1$, $D^{(j)}\Phi_x(0)
  = 0$.
\item[$(ii)$]  la fonction $\frac{1}{(r-k)!} D^{(r-k)}\Phi_x(0) z =
  \Psi_0(z)$ est un polyn{\^o}me homog{\`e}ne de degr{\'e} $r-k$ qui ne s'annule pas
  sur $\Omega$.
\end{itemize}
\hfill\hfill$\square$
\end{corollaire}

\begin{lemme}\label{lemme_non_nul_det}
Soit $x\in J$, et soit $K$ un compact contenu dans
$\Omega$.  Il existe une constante $\epsilon>0$ (d{\'e}pendant de $x$ et de $K$)
telle que pour tout $t$, $0<t\leq \epsilon$ et tout $\omega\in K$ on ait
$$ \det (x+t\omega) \not=0.$$
\end{lemme}
\proof
D'apr{\`e}s le lemme \ref{lemme_estim_det}, on peut {\'e}crire
$$ \det(x+t\omega) = 
t^{r-k}\Psi_0(\omega) + \sum_{j=1}^k t^{r-k+j}\Psi_j(\omega)$$
o{\`u} les $\Psi_j$ sont des polyn{\^o}mes homog{\`e}nes de degr{\'e} $ r-k+j$. La
fonction $\Psi_0$ garde un signe constant sur $\Omega$ (en fait le
signe de ${\det}_c(x)$ o{\`u} $c$ est l'idempotent donn{\'e} par le lemme \ref{lemme_unicite_c}). On va supposer qu'on est dans le cas o{\`u} ce
signe est $+$, l'autre cas se traitant de la m{\^e}me mani{\`e}re. Soit
$ c_0>0$ le minimum de $\Psi_0$ sur le compact $K$, et soit $C$ une
constante positive telle que pour tout $j$, $1\leq j\leq k$ on ait
$$\forall \omega\in K,\quad   |\Psi_j(\omega)|\leq C\ .$$
On a alors  pour $\omega\in K$ et $t>0$ 
$$ \det(x+t\omega)> t^{r-k}(c_0 -C\sum_{j=1}^k t^j)\ .$$
Posons $\displaystyle{\epsilon= \inf (1, \frac{1}{2}\frac{c_0}{kC})}$. Alors pour
$0\leq t\leq \epsilon$, on a 
$$c_0-C \sum_{j=1}^k t^j \geq c_0 - C k \epsilon\geq \frac{1}{2}c_0>0.$$
D'o{\`u} pour $\omega\in K$ et $0<t\leq \epsilon$, $\det(x+t\omega)>0$.
\hfill\hfill$\square$

\begin{proposition}\label{estim_det_c_causle}
Soit $\gamma(t)$, $0\leq t\leq 1$ une courbe causale d'origine $x\in J$.
\begin{itemize}
\item[$(i)$]  Il existe $\epsilon>0$ tel que $\forall t$, $0<t\leq
  \epsilon$, $\det (\gamma(t))\neq 0$.
\item[$(ii)$] Si $x$ est de rang $k$, il existe une constante $a_k\not=0$ telle que
$$\det (\gamma(t)) = a_k t^{r-k}+O(t^{r-k+1}).$$ 
\end{itemize}
\end{proposition}
\proof
Soit $ C =\{\dot\gamma(t),\; 0\leq t\leq 1\}$. Alors
$C$ est un compact de $\Omega$, et son enveloppe convexe
$\Gamma(C)$ {\'e}galement. {\'E}crivant
$$ \gamma(t) =t \int _0^1\dot{\gamma}(st) ds$$
on voit que $ \gamma(t) = t\omega(t)$, o{\`u} $\omega(t)\in \Gamma(C)$.
La partie $(i)$ est donc cons{\'e}quence du lemme \ref{lemme_non_nul_det}.\\
La partie $(ii)$ est cons{\'e}quence du corollaire \ref{corollaire_D_det}, $(ii)$.
\hfill\hfill$\square$\\

\begin{definition}
Soient $x$ et $y$ deux {\'e}l{\'e}ments quelconques de $J$. La quantit{\'e}
\begin{equation*}
\mu(x,y)=r-\mathrm{rang}(x-y)
\end{equation*}
est appel{\'e} l'{\it indice de transversalit{\'e}} de la paire $(x,y)$.
\end{definition}

Les deux cas extr{\^e}mes sont le cas o{\`u} $x-y$ est inversible, auquel cas
l'indice de transversalit{\'e} vaut 0, et le cas o{\`u} $x=y$, auquel cas l'indice de
transversalit{\'e} vaut $r$.\\

On va maintenant montrer que l'indice de transversalit{\'e} caract{\'e}rise les
orbites de $G(T_\Omega)$ dans $J\times J$.

\begin{remarque}\label{remarque_plong_J}
{\rm
D'apr{\`e}s (\ref{plongement_de_J}), le groupe conforme $G(T_\Omega)$ agit de
mani{\`e}re rationnelle sur l'alg{\`e}bre de Jordan $J$ si on la
r{\'e}alise dans sa compactification conforme $G(T_\Omega)/LN^-$. On peut
alors montrer que $\{g\in G(T_\Omega) \tq g(0)\in J\}=N^+LN^-$ et que
le stabilisateur de $0$ dans $G(T_\Omega)$ est $LN^-$. Pour toutes ses consid{\'e}rations voir \cite{Koufany}.
}
\end{remarque}

\begin{theoreme}\label{Invariance_indice_trans_J}
L'indice de transversalit{\'e} sur $J$ est invariant sous l'action du groupe
conforme $G(T_\Omega)$. Plus pr{\'e}cis{\'e}ment, si $g$ est une transformation conforme
de $J$ qui est d{\'e}finie en $x$ et $y$, alors
\begin{equation}\label{inv_mu_J}
\mu(g(x),g(y))=\mu(x,y).
\end{equation}
\end{theoreme}
\proof Il est imm{\'e}diat que l'indice de transversalit{\'e} est
invariant par translation.  Comme  le groupe $L$ conserve le rang d'un
{\'e}l{\'e}ment (cons{\'e}quence du th{\'e}or{\`e}me \ref{theo_inertie_Sylvester}), il conserve l'indice de
transversalit{\'e}. Reste donc {\`a} voir l'invariance sous
l'action de l'inversion. Soit donc $x, y$ deux {\'e}l{\'e}ments inversibles de
$J$. On utilise alors l'identit{\'e} de Hua sous la forme
\begin{equation*}
P(x^{-1}-y^{-1})=P(x)^{-1}P(x-y)P(y)^{-1}.
\end{equation*}
On en d{\'e}duit que $P(x^{-1}-y^{-1})$ et $P(x-y)$ ont m{\^e}me rang, mais
cela est {\'e}quivalent {\`a} dire (voir \cite[Proposition IV.3.1]{Faraut-Koranyi}) que
$x^{-1}-y^{-1}$ et $x-y$ ont m{\^e}me rang. D'o{\`u} le r{\'e}sultat.
\hfill\hfill$\square$\\

Abordons maintenant la r{\'e}ciproque.
\begin{lemme}\label{LN-transi_J_k} Soient $x, x'\in J$ tels que $\mu(x,0)=\mu(x',0)$. Alors il
  existe $h\in LN^-$ tel que $h(x)=x'$.
\end{lemme}
\proof
Posons $\mu(x,0)=\mu(x',0)=k$. L'{\'e}l{\'e}ment $x$ est donc de rang $\ell=r-k$. Il
existe un rep{\`e}re de Jordan
$(c_j)_{1\leq j\leq r}$ et  des nombres r{\'e}els  non nuls
$(\lambda_j)_{1\leq j\leq \ell}$ tels que
$x=\sum_{j=1}^\ell \lambda_j c_j$. 
Pour $\epsilon
>0$ assez petit, l'{\'e}l{\'e}ment $$x_\epsilon=\sum_{j=1}^\ell\lambda_jc_j+\sum_{j=\ell+1}^r\epsilon
c_j$$
 est inversible. Son image par $\widetilde{t}_v$  o{\`u}
 $v=\sum_{j=1}^\ell\frac{\lambda_j-1}{\lambda_j}c_j$  est donn{\'e}e par 
$$\begin{array}{ll}
\widetilde{t}_v(x_\epsilon)&=(s\circ t_v\circ s)(x_\epsilon)\\
                              &=(x_\epsilon^{-1}+v)^{-1}\\
                              &=\sum_{j=1}^\ell
c_j+\sum_{j=\ell+1}^r\epsilon c_j.
\end{array}$$
(voir (\ref{inversion}) et (\ref{transmutation}) pour ces diff{\'e}rentes notations).
Comme $\widetilde{t}_v$ est continue au voisinage de tout {\'e}l{\'e}ment
inversible, $\widetilde{t}_v(x_\varepsilon)$ tends vers $\sum_{j=1}^\ell c_j$ lorsque
$\varepsilon$ tend vers $0$. Par prolongement par continuit{\'e} on a donc $$\widetilde{t}_v(x)=\sum_{j=1}^\ell c_j.$$
L'{\'e}l{\'e}ment $x'$ est aussi de rang $\ell$, et il existe un syst{\`e}me
de Jordan $(c'_j)_{1\leq r}$ et des nombres r{\'e}els non nuls $(\lambda'_j)_{1\leq
  j\leq \ell}$ tels que $x'=\sum_{j=1}^\ell \lambda'_jc'_j$. On montre de
m{\^e}me qu'il existe une transformation $\widetilde{t}_{v'}$ telle que 
$$\widetilde{t}_{v'}(x')=\sum_{j=1}^\ell c'_j.$$
D'apr{\`e}s \cite[Theorem IV.2.5]{Faraut-Koranyi}, il existe un automorphisme
$k\in Aut(J)$ tel que $k(c_j)=c'_j$ pour tout $j$, $1\leq j\leq
r$. L'{\'e}l{\'e}ment $h=\widetilde{t}_{v'}^{-1}\circ k\circ \widetilde{t}_{v}\in
LN^-$ v{\'e}rifie alors $h(x)=x'$.
\hfill\hfill$\square$

\begin{theoreme}\label{J_k-homogene}
Soient $x, x_0, y, y_0 \in J$. Si $\mu(x,x_0)=\mu(y,y_0)$, alors il existe un
{\'e}l{\'e}ment $g\in N^+LN^-$ tel que $y=g(x)$ et $y_0=g(x_0)$.
\end{theoreme}
\proof 
La condition $\mu(x,x_0)=\mu(y,y_0)$ {\'e}tant {\'e}quivalente {\`a}
$\mu(x-x_0,0)=\mu(y-y_0,0)$, les deux {\'e}l{\'e}ments $x-x_0$ et $y-y_0$ ont le m{\^e}me
rang. D'apr{\`e}s le lemme \ref{LN-transi_J_k}, il existe $h\in LN^-$ tel que
$y-y_0=h(x-x_0)$. L'{\'e}l{\'e}ment $g=t_{y_0}\circ h\circ t_{-x_0}\in N^+LN^-$  v{\'e}rifie
$g(x)=y$ et $g(x_0)=y_0$ puisque $h$ fixe $0$ d'apr{\`e}s la remarque
\ref{remarque_plong_J}.
\hfill\hfill$\square$\\

Rappelons la d{\'e}finition de la transversalit{\'e} introduite dans \cite{Clerc-Orsted-1}.
On dira que deux {\'e}l{\'e}ments $\sigma$ et $\tau$ de  $S$ sont {\it transverses}  si,
\begin{equation*}
\sigma\top\tau \stackrel{def}{\Longleftrightarrow} \det(\sigma-\tau)\not=0.
\end{equation*}

Pour {\'e}tendre la d{\'e}finition de l'indice de transversalit{\'e} {\`a} $S\times
S$ nous avons besoin de la proposition suivante.

\begin{proposition}
Soient $\sigma, \tau\in S$, alors il existe un {\'e}l{\'e}ment $u$ de
$U$ tel que  $u(\sigma) $ et $u(\tau) $ soient transverses {\`a} $e$. De plus
l'entier 
$$\mu\bigl(c(u(\sigma) ),c(u(\tau)) \bigr)$$ ne d{\'e}pend pas de l'{\'e}l{\'e}ment
$u\in U$ choisi.
\end{proposition}
\proof
Soient $\sigma, \tau\in S$ et soit $\varsigma\in S$ tel que
$\varsigma\top \sigma$ et $\varsigma\top \tau$. On peut trouver $u\in U$, v{\'e}rifiant $u(\varsigma) =e$. Les
deux {\'e}l{\'e}ments  
$u(\sigma) $ et $u(\tau) $ sont alors transverses {\`a}
$e$ et par cons{\'e}quent appartiennent {\`a} $p(J)$.\\
Supposons
maintenant qu'il existe deux {\'e}l{\'e}ments $u_1, u_2\in U$ tels que $u_i(\sigma)$ et
$u_i(\tau) $ soient transverses {\`a} $e$ pour $i=1, 2$. Posons $c(u_i(\sigma))=x_i$
et $c(u_i(\tau))=y_i$, pour $i=1,2$. Alors $x_2=(cu_2u_1^{-1}c^{-1})(x_1)$ et
$y_2=(cu_2u_1^{-1}c^{-1})(y_1)$.
Comme $cu_2u_1^{-1}c^{-1}\in G(T_\Omega)$, on a d'apr{\`e}s (\ref{inv_mu_J}),
$$\begin{array}{ll}
\mu(x_1,y_1) &=\mu\bigl((cu_2u_1^{-1}c^{-1})(x_1),(cu_2u_1^{-1}c^{-1})(y_1)\bigr),\\
                &=\mu(x_2,y_2). 
\end{array} $$
\hfill\hfill$\square$

Ceci justifie la d{\'e}finition suivante 
\begin{definition} 
Soit $(\sigma,\tau)\in S\times S$. On appelle {\it indice de
transversalit{\'e}} de la paire $(\sigma,\tau)$ l'indice de transversalit{\'e} de
la paire 
$(c(u(\sigma)),c(u(\tau)))\in J\times J$ o{\`u} $u\in U$ est choisi 
tel que $u(\sigma)\top e$ et $u(\tau)\top e$. On note
\begin{equation*}
\mu(\sigma,\tau)=\mu\bigl(c(u(\sigma)),c(u(\tau))\bigr).
\end{equation*}
\end{definition}
En particulier, $\sigma$ et $\tau$ sont transverses si et seulement si
$\mu(\sigma,\tau)=0$.

\begin{remarque}{\rm
On peut montrer que $$\mu(\sigma,\tau)=k \iff \text{corang}P(\sigma-\tau)=k+\frac{k(k-1)}{2}d.$$
}
\end{remarque}

\begin{proposition} L'indice de transversalit{\'e} sur $S$ 
  est invariant par le groupe $G$,
\begin{equation*} 
\mu(g(\sigma),g(\tau))=\mu(\sigma,\tau),\;\;\forall g\in G.
\end{equation*} 
\end{proposition}
\proof Cela r{\'e}sulte du th{\'e}or{\`e}me \ref{Invariance_indice_trans_J} et du fait que $G(T_\Omega)=c\circ G\circ c^{-1}$. 

\hfill\hfill$\square$

\begin{theoreme}\label{th_S_k_orbite} L'indice de transversalit{\'e}
  caract{\'e}rise les orbites de l'action de $G$ sur $S\times S$. Plus
  pr{\'e}cis{\'e}ment, si $(\sigma, \sigma_0), (\tau, \tau_0)$ sont deux 
  {\'e}l{\'e}ments de $S\times S$ et si $\mu(\sigma,\sigma_0)=\mu(\tau,\tau_0)$ alors il
  existe $g\in G$ tel que $g(\sigma)=\tau$ et $g(\sigma_0)=\tau_0$.
\end{theoreme}
\proof Remarquons d'abord que d'apr{\`e}s l'invariance de $\mu$
on peut supposer
$\sigma_0=\tau_0=e$. Il suffit alors de montrer que si
$\mu(\sigma,e)=\mu(\tau,e)$ alors il existe $g\in G$ tel que $g(\sigma)=\tau$
et $g(e)=e$. Soit $\varsigma\in S$ tel que $\varsigma\top \sigma$,
$\varsigma\top \tau$ et $\varsigma\top e$. Soit $u\in U$ tel que
$u(\sigma)\in p(J)$, $u(\tau)\in p(J)$ et $u^{-1}(e)=\varsigma$. Posons
$x=c(u(\sigma))$, $y=c(u(\tau))$ et $z=c(u(e))$. La condition
$\mu(\sigma,e)=\mu(\tau,e)$ se traduit par $\mu(x,z)=\mu(y,z)$. La
proposition \ref{J_k-homogene} montre qu'il
existe $h\in N^+LN^-$ tel que $y=h(x)$ et $z=h(z)$. Donc l'{\'e}l{\'e}ment
$g=u^{-1}c^{-1}hcu\in G$ v{\'e}rifie $g(\sigma)=\tau$ et
$g(e)=e$. D'o{\`u} le th{\'e}or{\`e}me. \hfill\hfill$\square$

\begin{proposition}\label{estim_det_c_causale_S}
Soit $\sigma, \tau\in S$ et $[0,1]\ni t\mapsto \tau(t)\in S$ une courbe
causale d'origine $\tau$. Il existe $\epsilon>0$ tel que $\forall t$,
$0<t\leq \epsilon$, $\tau(t)$ est transverse {\`a} $\sigma$.
\end{proposition}
\proof Le r{\'e}sultat se d{\'e}duit par transformation de Cayley
de la proposition \ref{estim_det_c_causle}.
\hfill\hfill$\square$

\section{L'indice de Souriau : cas transverse}
Des r{\'e}sultats pr{\'e}c{\'e}dents, on d{\'e}duit que le groupe $G$ pr{\'e}serve la transversalit{\'e}, et agit transitivement sur 
\begin{equation*}
S^2_\top=\{(\sigma,\tau)\in S^2 \tq \sigma\top\tau\}.
\end{equation*}
Soit $S_\top(-e)$ l'ouvert d{\'e}fini par 
\begin{equation*}
S_\top(-e)= \{ \sigma\in S\tq \sigma\top (-e)\}.
\end{equation*}
Si $\sigma\in S_\top(-e)$, alors $ (\sigma-s e)$ est inversible pour
$s=-1$, et donc pour tout $s\in (-\infty,0]$. L'int{\'e}grale ({\`a} valeurs dans $\mathbb{J}$) 
\begin{equation*}
\log \sigma=\int_{-\infty}^0\big((se-\sigma)^{-1}-(s-1)^{-1}e\big)ds
\end{equation*}
est normalement convergente, et d{\'e}finit une fonction r{\'e}guli{\`e}re sur l'ouvert $S_\top(-e)$.
\begin{proposition}
Pour tout $\sigma\in S_\top(-e)$, on a
\begin{enumerate}
\item[$(i)$] $\exp(\log \sigma )=\sigma.$
\item[$(ii)$] $e^{\tr(\log \sigma )}=\det(\sigma).$
\item[$(iii)$] $\log \sigma^{-1} =-\log \sigma $.
\item[$(iv)$] $\log(k\sigma)=k\log \sigma $, pour tout $k\in Aut(J)$.\\
\end{enumerate}
\end{proposition}
\proof
Soit $\sigma\in S_\top(-e)$. Il existe une d{\'e}composition
de Peirce $e= \sum_{j=1}^r c_j$ de l'identit{\'e} adapt{\'e}e {\`a} $\sigma$,
c'est-{\`a}-dire telle que 
$$ \sigma = \sum_{j=1}^r \sigma_jc_j,\quad |\sigma_j|= 1, \ 1\leq j\leq
r\ .$$ 
Par hypoth{\`e}se, $\sigma_j\neq -1$ pour tout $j$ tel que $1\leq j\leq r$, et pour $s\in (-\infty,0]$
$$ (se-\sigma)^{-1}= \sum_{j=1}^r(s-\sigma_j)^{-1}c_j\ .$$
Par suite,
$$ \log \sigma = \sum_{j=1}^r \Big(\int_{-\infty}^0
\big((s-\sigma_j)^{-1}-(s-1)^{-1}\big)ds\Big) c_j\  .$$
On est donc ramen{\'e} au calcul, pour $z\in \mathbb{C}\setminus (-\infty, 0]$ de
l'int{\'e}grale
\begin{equation*}
\int_{-\infty}^0 \big((s-z)^{-1}-(s-1)^{-1}\big)ds.
\end{equation*}
Pour $z$ r{\'e}el strictement positif, un calcul {\'e}l{\'e}mentaire montre que cette
int{\'e}grale vaut $\log z$. Par prolongement analytique, on obtient que
l'int{\'e}grale co{\"\i}ncide avec la d{\'e}termination principale de $\log z$, pour
tout $z\in\mathbb{C}\setminus(-\infty,0]$. Donc, avec la m{\^e}me convention
\begin{equation}\label{log_1}
\log \sigma = \sum_{j=1}^r \log \sigma_j \ c_j,
\end{equation}
d'o{\`u} $(i)$ en calculant l'exponentielle des deux membres. Les autres
points de la proposition se d{\'e}duisent de (\ref{log_1}).
\hfill\hfill$\square$\\

Soit maintenant $\sigma$ et $\tau$ deux {\'e}l{\'e}ments transverses de $S$. Il
existe un {\'e}l{\'e}ment $u\in U$ tel que $ \tau = u(-e)$. L'{\'e}l{\'e}ment $u$ n'est pas
unique, mais tout autre {\'e}l{\'e}ment $u'$ satisfaisant $\tau =  u'(-e)$ peut
s'{\'e}crire $u'=uk$ o{\`u} $k\in U_{-e}=U_e=U\cap Aut(J)$. L'{\'e}l{\'e}ment
$ u^{-1}(\sigma)$ est transverse {\`a} $-e$, et donc on peut d{\'e}finir $\log
u^{-1}(\sigma)$. Pour  $u'=uk$, on obtient $ \log u'^{-1}(\sigma) = 
k^{-1}\log u^{-1}(\sigma)$ et donc $\tr(\log u^{-1}(\sigma))$ est
ind{\'e}pendant du repr{\'e}sentant $u$ choisi. En cons{\'e}quence, on pose 
\begin{equation*}
\Psi(\sigma, \tau) = \frac{1}{i}\tr(\log u^{-1}(\sigma)).
\end{equation*}
On d{\'e}finit ainsi une fonction sur $S^2_\top$ une fonction {\`a}
valeurs r{\'e}elles.
 De plus,
\begin{equation*}
e^{ \tr(\log u^{-1}(\sigma))}=\det(u^{-1}(\sigma)) = \chi(u)^{-1}\det(\sigma) = (-1)^r(\det \sigma) (\det \tau)^{-1}
\end{equation*}
de sorte que 
\begin{equation} \label{exp_Psi}
e^{2i\Psi(\sigma,\tau)}=(\det \sigma)^{2}(\det \tau) ^{-2}.
\end{equation}

\begin{proposition}
La fonction $\Psi : S^2_\top \to \mathbb{R}$ est r{\'e}guli{\`e}re. Soient $\sigma,
\tau\in S^2_\top$. Alors
\begin{enumerate}
\item[$(i)$] $\Psi(v(\sigma),v(\tau))=\Psi(\sigma,\tau)$, pour tout $v\in U$.
\item[$(ii)$]  $\Psi(\sigma,\tau)=-\Psi(\tau,\sigma)$.
\end{enumerate}
\end{proposition}
\proof Montrons la continuit{\'e} de $\Psi$ en un point
$(\sigma_0,\tau_0)\in S^2_\top$. Il est clair que l'{\'e}l{\'e}ment $u\in U$
tel que $\tau=u(-e)$ peut {\^e}tre choisi contin{\^u}ment pour $\tau$ proche
de $\tau_0$. Il suffit en effet d'utiliser le fait que $S$ est un espace
sym{\'e}trique sous l'action de $U$, et donc il y a une section locale
r{\'e}guli{\`e}re de l'application $U\to S$, $v\mapsto v^{-1}(\tau_0)$ qui vaut
$u$ en $-e$. D'apr{\`e}s le th{\'e}or{\`e}me de continuit{\'e} de l'int{\'e}grale
d'une fonction d{\'e}pendant de param{\`e}tres, il s'en suit que
$\Psi(\sigma,\tau)$ est bien une fonction continue de $(\sigma,\tau)$ au
voisinage de $(\sigma_0,\tau_0)$. La diff{\'e}rentiabilit{\'e} de $\Psi$ se
d{\'e}montre de m{\^e}me.\\
Soit $v\in U$. Alors, si $\tau=u(-e)$, on a $v(\tau)=(v\circ u)(-e)$, de
sorte que
\begin{equation*}
\Psi(v(\sigma),v(\tau))=\frac{1}{i}\tr\Bigl[\log\bigl((v\circ
u)^{-1}(v(\sigma))\bigr)\Bigr]=\frac{1}{i}\tr\bigl(\log u^{-1}(\sigma)\bigr)=\Psi(\sigma,\tau).
\end{equation*}
Pour d{\'e}montrer $(ii)$, on peut compte-tenu de $(i)$ supposer que
$\tau=-e$. Choisissons une racine carr{\'e}e $\sigma^{\frac{1}{2}}$ de
$\sigma$. Comme $i\sigma^{\frac{1}{2}}\in S$, l'{\'e}l{\'e}ment
$P(i\sigma^{\frac{1}{2}})$ appartient {\`a} $U$, et satisfait
$P(i\sigma^{\frac{1}{2}})(-e)=\sigma$. Par suite
\begin{equation*}
\Psi(-e,\sigma)=\frac{1}{i}\tr\Bigl[\log\bigl(P(i\sigma^{\frac{1}{2}})^{-1}(-e)\bigr)\Bigr]=\frac{1}{i}\tr\bigl(\log\sigma^{-1}\bigr)=-\frac{1}{i}\tr(\log\sigma)=-\Psi(\sigma,-e)
\end{equation*}
o{\`u} on a utilis{\'e} le fait que $\sigma^{-1}=\overline{\sigma}$ pour tout
{\'e}l{\'e}ment de $S$. D'o{\`u} $(ii)$.
\hfill\hfill$\square$

Pour  $\theta\in\mathbb{R}$, $\theta\not\in \pi+2\pi\mathbb{Z}$, on note
$\{\theta\}$ le r{\'e}el caract{\'e}ris{\'e} par
\begin{equation*}
\{\theta\}\equiv \theta \;[2\pi],\quad -\pi<\{\theta\}<\pi
\end{equation*}
(une version de la {\it d{\'e}termination principale} de $\theta$).

\begin{proposition}
Soit $(c_j)_{1\leq j\leq r}$ une rep{\`e}re de Jordan de $J$, et soit
\begin{equation*}
\sigma=\sum_{j=1}^re^{i\theta_j}c_j,\;\;\; \tau=\sum_{j=1}^re^{i\phi_j}c_j
\end{equation*}
deux {\'e}l{\'e}ments transverses de $S$ (c'est-{\`a}-dire tels que
$\theta_j\not\equiv \phi_j\;[2\pi]$ pour tout $j$, $1\leq j\leq
r$). Alors
\begin{equation}\label{equation_Psi}
\Psi(\sigma,\tau)=\sum_{j=1}^r\{\theta_j-\phi_j+\pi\}.
\end{equation}
\end{proposition}
\proof
Posons $u=\exp\bigl[iL(\sum_{j=1}^r(\phi_j-\pi)c_j)\bigr]\in U$. On a
$u(-e)=\tau$ et
$u^{-1}(\sigma)=\sum_{j=1}^re^{i(\theta_j-\phi_j+\pi)}c_j$. D'o{\`u} le
r{\'e}sultat, d'apr{\`e}s (\ref{log_1}).
\hfill\hfill$\square$\\

On {\'e}tend la relation de transversalit{\'e} {\`a}  $\widetilde{S}$, en posant, pour
deux {\'e}l{\'e}ments 
$(\sigma,\theta), (\tau,\phi)$ de  $\widetilde{S}$
\begin{equation*}
(\sigma,\theta)\top (\tau,\phi)\Longleftrightarrow \sigma\top
\tau
\end{equation*}
et on note $\widetilde{S}^2_\top=\{(\widetilde{\sigma},\widetilde{\tau})\in\widetilde{S}^2
\tq \widetilde{\sigma}\top \widetilde{\tau}\}$.
Pour un couple d'{\'e}l{\'e}ments transverses $ \widetilde{\sigma}= (\sigma,\theta),\widetilde{\tau}=
(\tau,\phi)$, on d{\'e}finit l'{\it indice de Souriau}
$m(\widetilde{\sigma},\widetilde{\tau})$ par la formule
\begin{equation}\label{indice_S_L}
m(\widetilde \sigma,\widetilde \tau)  = \frac{1}{\pi}\big[\Psi(\sigma,\tau)-r(\theta-\phi)\big] .
\end{equation}
On remarque tout de suite que $m$ est antisym{\'e}trique,
\begin{equation*}
m(\widetilde{\sigma},\widetilde{\tau})+m(\widetilde{\tau},\widetilde{\sigma})=0.
\end{equation*}
De plus. d'apr{\`e}s (\ref{exp_Psi})
\begin{equation*}
e^{2i\pi m(\widetilde{\sigma},\widetilde{\tau})}=
e^{-2ir\theta}e^{2ir\phi}(\det\,\sigma)^{2}(\det \tau)^{-2}=1,
\end{equation*}
de sorte que $ m(\widetilde{\sigma},\widetilde{\tau})\in \mathbb{Z}$.\\
 Enfin la fonction $m$ est une fonction continue (donc localement constante) sur
$\widetilde{S}^2_\top$, d'apr{\`e}s la continuit{\'e} de la fonction $\Psi$.

\begin{proposition}\label{invarince_m}
La fonction $ m$ est invariante sous l'action
du groupe $\Gamma$.
\end{proposition}
\proof Soit $\mathcal{U}$ un voisinage assez petit et simplement
connexe de l'{\'e}l{\'e}ment neutre  de $G$. On peut trouver une d{\'e}termination $\varphi(g,z)$ de
l'argument de
$j(g,z)$ qui soit d{\'e}finie sur $\mathcal{U}\times \overline{\mathcal{D}}$ et
d{\'e}pendant contin{\^u}ment du couple $(g,z) $. Soient $\widetilde{\sigma} = (\sigma,\theta), \widetilde{\tau} = (\tau,
\phi)$ deux {\'e}l{\'e}ments de $\widetilde{S}$, et supposons $ \sigma\top \tau$.  Soit
$g\in\mathcal{U}$, alors
$g(\sigma)\top g(\tau)$, et on a
$$
m\bigl((g,\varphi(g,\cdot))(\sigma,\theta),(g,\varphi(g,\cdot))(\tau,\phi)\bigr)
= \frac{1}{\pi}[\Psi(g(\sigma), g(\tau)) - r(\theta-\phi)
-(\varphi(g,\sigma)-\varphi(g,\tau))]$$
 ce qui montre que $g\longmapsto
m\bigl((g,\varphi(g,\cdot))(\sigma,\theta),(g,\varphi(g,\cdot))(\tau,\phi)\bigr)$
est une fonction continue. Comme elle est {\`a} valeurs dans $\mathbb{Z}$
elle est donc constante sur $\mathcal{U}$. Par ailleurs on v{\'e}rifie
facilement, d'apr{\`e}s (\ref{action_par_T}), que $m(T\cdot\widetilde{\sigma}, T\cdot\widetilde{\tau}) = m(\widetilde{\sigma},\widetilde{\tau})$.
Le voisinage $\{(g,\varphi(g,\cdot) \tq g\in\mathcal{U}\}$ du neutre de $\Gamma$  et l'{\'e}l{\'e}ment $T$ engendrent le groupe
$\Gamma$, la fonction $m$ est bien invariante par $\Gamma$.
\hfill\hfill$\square$\\

Soit maintenant
\begin{equation*}
S^3_\top=\{(\sigma_1,\sigma_2,\sigma_3)\in S^3 \tq \sigma_j\top \sigma_k,\;
1\leq j\not=k\leq 3\}.
\end{equation*}
Fixons un rep{\`e}re de Jordan $(c_j)_{1\leq j\leq r}$ et posons pour
tout entier $k$ tel que $0\leq k\leq r$
\begin{equation*}
\varepsilon_0=-e,\; \varepsilon_k=\sum_{j=1}^kc_j-\sum_{j=k+1}^rc_j,\; \varepsilon_r=e.
\end{equation*}
Le groupe $G$ agit sur $S^3_\top$ et d'apr{\`e}s \cite[Theoreme
4.3]{Clerc-Orsted-1} cette action admet exactement $r+1$ orbites
repr{\'e}sent{\'e}es par la famille $(e,-e,-i\varepsilon_k)$, $0\leq k\leq r$. \\
Soit $(\sigma_1,\sigma_2,\sigma_3)\in S^3_\top$. L'{\it indice de Maslov} $\imath(\sigma_1,\sigma_2,\sigma_3)$
 est  d{\'e}fini dans \cite[(21)]{Clerc-Orsted-1} par
\begin{equation*}
\imath(\sigma_1,\sigma_2,\sigma_3)=2k-r
\end{equation*}
o{\`u} $k$ est l'unique entier tel que $(\sigma_1,\sigma_2,\sigma_3)$ soit
conjugu{\'e} {\`a} $(e,-e,-i\varepsilon_k)$.\\

Posons
\begin{equation*}
\widetilde{S}^3_\top=\{(\widetilde{\sigma}_1,\widetilde{\sigma}_2,\widetilde{\sigma}_3)\in S^3 \tq \widetilde{\sigma}_j\top \widetilde{\sigma}_k,\;
1\leq j\not=k\leq 3\}.
\end{equation*}
Le th{\'e}or{\`e}me suivant montre que l'indice de Souriau est en fait une
primitive de l'indice de Maslov.

\begin{theoreme}[Formule de Leray]\label{Formule_Leray} Soit
  $(\widetilde{\sigma}_1,\widetilde{\sigma}_2,\widetilde{\sigma}_3)\in
      \widetilde{S}^3_\top$. Alors la quantit{\'e}
\begin{equation*}
m(\widetilde{\sigma}_1,\widetilde{\sigma}_2)+m(\widetilde{\sigma}_2,\widetilde{\sigma}_3)+m(\widetilde{\sigma}_3,\widetilde{\sigma}_1)
\end{equation*}
ne d{\'e}pend que des projections $\sigma_1,\sigma_2,\sigma_3$ et est {\'e}gale
{\`a} l'indice de Maslov $\imath(\sigma_1,\sigma_2,\sigma_3)$. Autrement dit
\begin{equation*}
m(\widetilde{\sigma}_1,\widetilde{\sigma}_2)+m(\widetilde{\sigma}_2,\widetilde{\sigma}_3)+m(\widetilde{\sigma}_3,\widetilde{\sigma}_1)=\imath(\sigma_1,\sigma_2,\sigma_3).
\end{equation*}

\end{theoreme}
\proof
Par un calcul direct, en utilisant (\ref{indice_S_L}), on montre que
$$m(\widetilde{\sigma}_1,\widetilde{\sigma}_2)+m(\widetilde{\sigma}_2,\widetilde{\sigma}_3)+m(\widetilde{\sigma}_3,\widetilde{\sigma}_1)=\Psi(\sigma_1,\sigma_2)+\Psi(\sigma_2,\sigma_3)+\Psi(\sigma_3,\sigma_1).$$
Cette somme est donc un
entier qui ne d{\'e}pendent que des trois projections correspondantes. Soit
$\jmath(\sigma_1,\sigma_2,\sigma_3)$ l'entier ainsi d{\'e}fini. C'est une
fonction continue (donc en fait localement constante) de
$S^3_\top$. Le groupe $G$ est connexe et op{\`e}re contin{\^u}ment sur $S^3_\top$. Cela implique que la fonction $\jmath$ est invariante sous l'action de
$G$. Pour montrer que $\jmath$ et $\imath$ co{\"\i}ncident, il suffit donc de
les comparer sur chaque classe de conjugaison de triplets sous l'action de
$G$. Soit donc $k$ un entier tel que $0\leq k\leq r$. On doit calculer
$\jmath(e,-e,-i\varepsilon_k)$. On a
\begin{equation*}
e=\sum_{j=1}^re^{i0}c_j,\;\; -e=\sum_{j=1}^re^{i\pi}c_j,\;\; -i\varepsilon_k=\sum_{j=1}^ke^{-i\frac{\pi}{2}}c_j+\sum_{j=k+1}^re^{i\frac{\pi}{2}}c_j,
\end{equation*}
qu'on rel{\`e}ve en les points
\begin{equation*}
\widetilde{e}=(e,0),\;\; \widetilde{-e}=(-e,\pi),\;\; \widetilde{-i\varepsilon}_k=(-i\varepsilon_k,\frac{(r-2k)\pi}{2r}).
\end{equation*}
Un calcul simple en utilisant la formule (\ref{equation_Psi}) donne
\begin{equation*}
m(\widetilde{e},\widetilde{-e})=r,\;\;
m(\widetilde{-e},\widetilde{-i\varepsilon}_k)=-r,\;\;
m(\widetilde{-i\varepsilon}_k,\widetilde{e})= 2k-r.
\end{equation*}
D'o{\`u} 
\begin{equation*}
 m(\widetilde{\sigma}_1,\widetilde{\sigma}_2)+m(\widetilde{\sigma}_2,\widetilde{\sigma}_3)+m(\widetilde{\sigma}_3,\widetilde{\sigma}_1)=2k-r=\imath(\sigma_1,\sigma_2,\sigma_3).
\end{equation*}
\hfill\hfill$\square$

\section{L'extension de l'indice de Souriau au cas non transverse}

Dans \cite{Clerc2} il est montr{\'e} que l'indice de Maslov peut {\^e}tre
d{\'e}fini pour tout triplet $(\sigma_1, \sigma_2, \sigma_3)$ de $S\times
S\times S$. Cet indice, encore not{\'e} $\imath(\sigma_1, \sigma_2,
\sigma_3)$ prolonge l'indice d{\'e}fini sur $S^3_\top$ et poss{\`e}de encore
les propri{\'e}t{\'e}s principales, {\`a} savoir l'invariance par le groupe $G$,
l'antisym{\'e}trie et surtout la relation de cocycle que nous utiliserons
dans la suite sous la forme 
\begin{equation}\label{rel_cocy_g}
\imath(\sigma_1,\sigma_2,\sigma_3)-\imath(\sigma_1,\sigma_2,\sigma_4)+\imath(\sigma_1,\sigma_3,\sigma_4)-\imath(\sigma_2,\sigma_3,\sigma_4)=0
\end{equation}
pour tous $\sigma_1, \sigma_2, \sigma_3, \sigma_4\in S$.\\
Il est donc naturel de chercher {\`a} {\'e}tendre la d{\'e}finition de l'indice
de Souriau aux couples non transverses pour obtenir une primitive de
l'indice de Maslov {\'e}tendu. Pour l'extension de la primitive $m$ aux couples
non transverses, on va utiliser une id{\'e}e due {\`a} de Gosson
\cite{deGosson}.\\ 
Montrons d'abord l'unicit{\'e} de l'extension. Supposons que la fonction $m$ se
prolonge {\`a} $\widetilde{S}\times \widetilde{S}$ en une primitive
antisym{\'e}trique not{\'e}e encore $m$ du cocycle de Maslov, c'est-{\`a}-dire telle
que l'on ait
pour tout $(\widetilde{\sigma}_1,\widetilde{\sigma}_2,\widetilde{\sigma}_3)\in \widetilde{S}\times \widetilde{S}\times \widetilde{S}$,
\begin{equation*}
\imath(\sigma_1,\sigma_2,\sigma_3)=
m(\widetilde{\sigma}_1,\widetilde{\sigma}_2)+m(\widetilde{\sigma}_2,\widetilde{\sigma}_3)+m(\widetilde{\sigma}_3,\widetilde{\sigma}_1).
\end{equation*}
On en d{\'e}duit que n{\'e}cessairement
\begin{equation}\label{pre_def_ext_m}
m(\widetilde{\sigma}_1,\widetilde{\sigma}_2)=\imath(\sigma_1,\sigma_2,\sigma_3)
+m(\widetilde{\sigma}_1,\widetilde{\sigma}_3)
+m(\widetilde{\sigma}_3,\widetilde{\sigma}_2).
\end{equation}
{\'E}tant donn{\'e}s deux {\'e}l{\'e}ments $\widetilde{\sigma}_1,
\widetilde{\sigma}_2\in \widetilde{S}$, on peut toujours trouver
$\widetilde{\sigma}_3\in \widetilde{S}$ tels que $\widetilde{\sigma}_3\top
\widetilde{\sigma}_1$ et $\widetilde{\sigma}_3\top
\widetilde{\sigma}_2$. Le membre de droite de l'{\'e}galit{\'e}
(\ref{pre_def_ext_m}) est bien d{\'e}fini et il y a donc unicit{\'e} de l'extension possible.

\begin{theoremedefinition}\label{Def_ind_S_L}
Soient $\widetilde{\sigma}_1,\;
\widetilde{\sigma}_2$ deux {\'e}l{\'e}ments de
$\widetilde{S}$ et soit $\widetilde{\tau}$ un {\'e}l{\'e}ment quelconque de
$\widetilde{S}$ tel que $\widetilde{\tau}\top \widetilde{\sigma}_1$ et
$\widetilde{\tau}\top \widetilde{\sigma}_2$. Soient $\sigma_1,\sigma_2,\tau$ les projections respectives de
$\widetilde{\sigma}_1, \widetilde{\sigma}_2, \widetilde{\tau}$.
L'entier 
\begin{equation*}
\imath(\sigma_1,\sigma_2,\tau)+m(\widetilde{\sigma}_1,\widetilde{\tau})+m(\widetilde{\tau},\widetilde{\sigma}_2)
\end{equation*} ne d{\'e}pends pas du choix de $\widetilde{\tau}$ et est appel{\'e}
l'{\it indice de Souriau} du couple
$(\widetilde{\sigma}_1,\widetilde{\sigma}_2)$. Lorsque
$\sigma_1\top\sigma_2$, il co{\"\i}ncide avec $m(\widetilde{\sigma}_1,\widetilde{\sigma}_2)$.
\end{theoremedefinition}
\proof Soient
$\widetilde{\sigma}_1,\widetilde{\sigma}_2$ deux {\'e}l{\'e}ments quelconques de 
$\widetilde{S}$. Alors on peut trouver deux {\'e}l{\'e}ments
$\widetilde{\tau}_1,\widetilde{\tau}_2$ tels que
\begin{equation*}
\widetilde{\tau}_1\top \widetilde{\sigma}_1,\;
 \widetilde{\tau}_1\top\widetilde{\sigma}_2,\;\; 
\widetilde{\tau}_2\top \widetilde{\sigma}_1, 
\; \widetilde{\tau}_2\top \widetilde{\sigma}_2,\;\;
\mathrm{et}\;\;
 \widetilde{\tau}_1\top \widetilde{\tau}_2.
\end{equation*}
D'apr{\`e}s le th{\'e}or{\`e}me \ref{Formule_Leray}
appliqu{\'e} aux triplets transverses 
$(\widetilde{\sigma}_1,\widetilde{\tau}_1, \widetilde{\tau}_2)$ et
$(\widetilde{\sigma}_2,\widetilde{\tau}_1, \widetilde{\tau}_2)$, on obtient
\begin{eqnarray*}\label{m_g_1}
\imath(\sigma_1,\tau_1,\tau_2)=m(\widetilde{\sigma}_1,\widetilde{\tau}_1)+m(\widetilde{\tau}_1,\widetilde{\tau}_2)+m(\widetilde{\tau}_2,\widetilde{\sigma}_1),\\
\imath(\sigma_2,\tau_1,\tau_2)=m(\widetilde{\sigma}_2,\widetilde{\tau}_1)+m(\widetilde{\tau}_1,\widetilde{\tau}_2)+m(\widetilde{\tau}_2,\widetilde{\sigma}_2).
\end{eqnarray*}
D'o{\`u}
\begin{equation*}\label{m_g_2}
\imath(\sigma_1,\tau_1,\tau_2)-\imath(\sigma_2,\tau_1,\tau_2)=
m(\widetilde{\sigma}_1,\widetilde{\tau}_1)+m(\widetilde{\tau}_2,\widetilde{\sigma}_1)
-m(\widetilde{\sigma}_2,\widetilde{\tau}_1)-m(\widetilde{\tau}_2,\widetilde{\sigma}_2).
\end{equation*}
En utilisant la relation de cocycle (\ref{rel_cocy_g}), cette derni{\`e}re
{\'e}galit{\'e} devient
\begin{equation*}
\imath(\sigma_1,\sigma_2,\tau_2)+m(\widetilde{\sigma}_1,\widetilde{\tau}_2)+m(\widetilde{\tau}_2,\widetilde{\sigma}_2)
=
\imath(\sigma_1,\sigma_2,\tau_1)+m(\widetilde{\sigma}_1,\widetilde{\tau}_1)+m(\widetilde{\tau}_1,\widetilde{\sigma}_2).
\end{equation*}
D'o{\`u} la premi{\`e}re assertion du th{\'e}or{\`e}me.\\
Lorsque $\sigma_1\top\sigma_2$, la formule de Leray (th{\'e}or{\`e}me
\ref{Formule_Leray}) appliqu{\'e}e au triplet
$(\widetilde{\sigma}_1, \widetilde{\sigma}_2, \widetilde{\tau})$ montre que
$m(\widetilde{\sigma}_1,\widetilde{\sigma}_2)=\imath(\sigma_1,\sigma_2,\tau)+m(\widetilde{\sigma}_1,\widetilde{\tau})+m(\widetilde{\tau},\widetilde{\sigma}_2)$.
D'o{\`u} la deuxi{\`e}me assertion du th{\'e}or{\`e}me.

\hfill\hfill$\square$\\

Soient $\widetilde{\sigma}_1, \widetilde{\sigma}_2$ deux {\'e}l{\'e}ments de
$\widetilde{S}$ de projections respectives $\sigma_1$ et $\sigma_2$. On pose dans la suite
\begin{equation}\label{Formule_S-L-generale}
m(\widetilde{\sigma}_1,\widetilde{\sigma}_2)=\imath(\sigma_1,\sigma_2,\tau)+m(\widetilde{\sigma}_1,\widetilde{\tau})+m(\widetilde{\tau},\widetilde{\sigma}_2),
\end{equation}
o{\`u} $\widetilde{\tau}\in \widetilde{S}$ est un {\'e}l{\'e}ment dont la
projection $\tau$ est transverse {\`a} $\sigma_1$ et $\sigma_2$. Cette
quantit{\'e} est encore appel{\'e}e l'indice de Souriau du couple $(\widetilde{\sigma}_1, \widetilde{\sigma}_2)$.

\begin{remarque}{\rm On v{\'e}rifie facilement les assertions suivantes :
    \begin{enumerate}
    \item[$(i)$] l'extension $m$ est antisym{\'e}trique.
    \item[$(ii)$] l'extension $m$ est invariante sous l'action du groupe $\Gamma$.
    \end{enumerate}
}
\end{remarque}

\begin{theoreme}[Formule de Leray {\'e}tendue]
Soient $\widetilde{\sigma}_1$, $\widetilde{\sigma}_2,$
$\widetilde{\sigma}_3)\in \widetilde{S}$ de projections
respectives $\sigma_1,\sigma_2,\sigma_3$, alors
\begin{equation}\label{Formule_Leray_Gen}
m(\widetilde{\sigma}_1,\widetilde{\sigma}_2)+m(\widetilde{\sigma}_2,\widetilde{\sigma}_3)+m(\widetilde{\sigma}_3,\widetilde{\sigma}_1)=\imath(\sigma_1,\sigma_2,\sigma_3).
\end{equation}
\end{theoreme}
\proof
Choisissons un {\'e}l{\'e}ment $\widetilde{\tau}\in \widetilde{S}$ tel que
$\widetilde{\tau}\top \widetilde{\sigma}_1$, $\widetilde{\tau}\top
\widetilde{\sigma}_2$ et $\widetilde{\tau}\top
\widetilde{\sigma}_3$. Alors d'apr{\`e}s la d{\'e}finition (\ref{Formule_S-L-generale}) de l'indice $m$
{\'e}tendu on a
\begin{eqnarray*}
m(\widetilde{\sigma}_1,\widetilde{\sigma}_3)=\imath(\sigma_1,\sigma_3,\tau)+m(\widetilde{\sigma}_1,\widetilde{\tau})+m(\widetilde{\tau},
\widetilde{\sigma}_3),\\
m(\widetilde{\sigma}_2,\widetilde{\sigma}_3)=\imath(\sigma_2,\sigma_3,\tau)+m(\widetilde{\sigma}_2,\widetilde{\tau})+m(\widetilde{\tau},\widetilde{\sigma}_3).
\end{eqnarray*}
En utilisant encore une fois la relation de cocycle (\ref{rel_cocy_g}), on obtient 
\begin{equation*}
\imath(\sigma_1,\sigma_2,\tau)+m(\widetilde{\sigma}_1,\widetilde{\tau})+m(\widetilde{\tau},\widetilde{\sigma}_2)=
\imath(\sigma_1,\sigma_2,\sigma_3)+m(\widetilde{\sigma}_1,\widetilde{\sigma}_3)+m(\widetilde{\sigma}_3,\widetilde{\sigma}_2).
\end{equation*}
Le membre de gauche de cette derni{\`e}re {\'e}galit{\'e} {\'e}tant par
d{\'e}finition $m(\widetilde{\sigma}_1,\widetilde{\sigma}_2)$, le th{\'e}or{\`e}me
s'en d{\'e}duit en utilisant l'antisym{\'e}trie de l'indice de Souriau {\'e}tendu.
\hfill\hfill$\square$\\

On va maintenant donner l'expression ``en coordonn{\'e}es'' de l'indice de
Souriau. Nous aurons besoin d'une notation particuli{\`e}re, qui n'est en fait
qu'une version de l'indice de Maslov pour le cas du cercle unit{\'e}.\\

Soient donn{\'e}s $e^{i\theta}, e^{i\phi}$ deux nombres complexes de module 1
distincts. On notes $]e^{i\theta}, e^{i\phi}[$ l'arc de cercle
d'extr{\'e}mit{\'e}s $e^{i\theta}$ et $e^{i\phi}$ form{\'e} des points obtenus en
partant de $e^{i\theta}$ dans le sens trigonom{\'e}trique pour arriver {\`a}
$e^{i\phi}$.\\

Pour $e^{i\theta}$, $e^{i\varphi}$ et $e^{i\phi}$ trois points du cercle
orient{\'e}, on pose 
$$
\text{ord}(e^{i\theta}, e^{i\varphi}, e^{i\phi}) =
\left\lbrace
\begin{array}{rl}
0& \text{si deux des points}\; e^{i\theta}, e^{i\varphi}, e^{i\phi}\; \text{co{\"\i}ncident},\\
1& \text{si}\; e^{i\varphi}\in ] e^{i\theta}, e^{i\phi} [\; , \\
-1& \text{si}\; e^{i\varphi}\in ] e^{i\phi}, e^{i\theta} [\; .
\end{array} \right.
$$
\begin{proposition} Soient $\sigma_1=\sum_{j=1}^r e^{i\theta_j}c_j$,
  $\sigma_2=\sum_{j=1}^r e^{i\varphi_j}c_j$ et $\sigma_3=\sum_{j=1}^r
  e^{i\phi_j}c_j$ trois points de $S$. Alors
\begin{equation}\label{indice_i_en_coord}
\imath(\sigma_1,\sigma_2,\sigma_3)=\sum_{j=1}^r\mathrm{ord}(e^{i\theta_j}, e^{i\varphi_j}, e^{i\phi_j}).
\end{equation}
\end{proposition}
\proof Cette formule est implicitement contenue dans \cite{Clerc2}.
\hfill\hfill$\square$\\

Fixons maintenant un rep{\`e}re de Jordan $(c_j)_{1\leq j\leq r}$. 
\begin{proposition}
Soient
$$\widetilde{\sigma}_1=(\sum_{j=1}^\ell e^{i\theta_j}c_j+\sum_{j=\ell+1}^r e^{i\theta_j}
c_j,\theta),\;\;
\widetilde{\sigma}_2=(\sum_{j=1}^\ell e^{i\theta_j}c_j+\sum_{j=\ell+1}^r e^{i\varphi_j}
c_j,\varphi), 
$$ deux {\'e}l{\'e}ments tels que $\theta_j-\varphi_j\notin2\pi\mathbb{Z}$ pour
$\ell+1\leq j\leq r$. Alors
\begin{equation}
  \label{eq:primitive}
 m(\widetilde{\sigma}_1,\widetilde{\sigma}_2)=\frac{1}{\pi}
\Bigl[\sum_{j=\ell+1}^r\{\theta_j-\varphi_j+\pi\}-r(\theta-\varphi)\Bigr]. 
\end{equation}
\end{proposition}
\proof
Pour un r{\'e}el assez petit $\phi>0$, l'{\'e}l{\'e}ment
$$\widetilde{\sigma}_3=(\sum_{j=1}^\ell
e^{i(\theta_j+\phi)}c_j+\sum_{j=\ell+1}^r e^{i(\varphi_j+\phi)}
c_j,\varphi+\phi),
$$
est transverse {\`a} $\widetilde{\sigma}_1$ et {\`a} $\widetilde{\sigma}_2$. D'une
part, en appliquant la formule (\ref{indice_i_en_coord}), l'indice de Maslov des projections $\sigma_1$, $\sigma_2$ et $\sigma_3$ est
$$\imath(\sigma_1,\sigma_2,\sigma_3)=r-\ell.$$
D'autre part, en utilisant (\ref{equation_Psi}) et (\ref{indice_S_L}), on obtient
$$
\begin{array}{ll}
m(\widetilde{\sigma}_1,\widetilde{\sigma}_3)&=\frac{1}{\pi}
\Bigl[\Bigl(\sum_{j=1}^\ell\{-\phi+\pi\}+\sum_{j=\ell+1}^r\{\theta_j-\varphi_j-\phi+\pi\}\Bigr)-r(\theta-(\varphi+\phi))\Bigr],\\
&=\ell+\frac{1}{\pi}\Bigl[\sum_{j=\ell+1}^r\{\theta_j-\varphi_j+\pi\}-r(\theta-\varphi)\Bigr],\\
\text{et}\\
m(\widetilde{\sigma}_3,\widetilde{\sigma}_2)&=\frac{1}{\pi}
\Bigl[\Bigl(\sum_{j=1}^\ell\{\phi+\pi\}+\sum_{j=\ell+1}^r\{\phi+\pi\}\Bigr)-r\phi\Bigr],\\
&=-r.
\end{array}
$$
Finalement, d'apr{\`e}s (\ref{Formule_S-L-generale}), l'indice de Souriau
$m(\widetilde{\sigma}_1,\widetilde{\sigma}_2)$ est donn{\'e} par
\begin{equation}
 m(\widetilde{\sigma}_1,\widetilde{\sigma}_2)=\frac{1}{\pi}
\Bigl[\sum_{j=\ell+1}^r\{\theta_j-\varphi_j+\pi\}-r(\theta-\varphi)\Bigr]. 
\end{equation}
\hfill\hfill$\square$\\

\noindent En particulier si $\widetilde{\sigma}_1=(-e,-\pi)$ et
$\widetilde{\sigma}_2=(-\sum_{j=1}^\ell c_j+\sum_{j=\ell+1}^re^{i\varphi_j}c_j,\varphi)$,
o{\`u} 
\begin{itemize}
\item[$(i)$] $-\pi<\varphi_j<\pi$, $\forall j$, $\ell+1\leq j\leq r$, et
\item[$(ii)$] $r\varphi=-\ell\pi+\sum_{j=\ell+1}^r\varphi_j+2k\pi$, avec
  $k\in\mathbb{Z}$,
\end{itemize}
alors
\begin{equation}
  \label{eq:m_en_coord_simple}
  m(\widetilde{\sigma}_1,\widetilde{\sigma}_2)=2k+r-\ell=2k+r-\mu(\sigma_1,\sigma_2).\\
\end{equation}

\begin{remarque}{\rm
Pour $\widetilde{\sigma}, \widetilde{\tau}\in \widetilde{S}$, on a
 \begin{equation}\label{action_T_m}
 m(\widetilde{\sigma},T\cdot\widetilde{\tau})=-m(T\cdot\widetilde{\sigma},\widetilde{\tau})=m(\widetilde{\sigma},\widetilde{\tau})+2
 \end{equation}
}
\end{remarque}

\section{Le cycle de Maslov}
Soit $\sigma_0\in S$. On note $\Sigma(\sigma_0)$ l'ensemble des
{\'e}l{\'e}ments de $S$ qui ne sont pas transverses {\`a} $\sigma_0$, c'est-{\`a}-dire
\begin{equation*}
\Sigma(\sigma_0)=S\setminus S_\top(\sigma_0)=\{\sigma\in S \tq \det(\sigma-\sigma_0)=0\}.
\end{equation*}
Cet ensemble est appel{\'e} le {\it cycle de Maslov\/} associ{\'e} {\`a} $\sigma_0$.\\
On va, dans ce paragraphe {\'e}tudier sa structure de vari{\'e}t{\'e} stratifi{\'e}e
(au sens alg{\'e}brique).
\subsection{La vari{\'e}t{\'e} des {\'e}l{\'e}ments de $J$ de corang $k$}
Reprenons les notations du paragraphe 4.
\begin{proposition}\label{lemme_Frobenius}
Soit $c$ un idempotent, et $J=J(c,1)\oplus J(c,\frac{1}{2})\oplus J(c,0)$
la d{\'e}composition de Peirce relative {\`a} cet idempotent. L'application donn{\'e}e
par 
\begin{equation*}
J(c,1)\times J(c,\frac{1}{2})\times J(c,0)\to J \; : \;
(x_1,x_{\frac{1}{2}},x_0)\mapsto \tau_c(x_{\frac{1}{2}})(x_1+x_0)
\end{equation*}
est un diff{\'e}omorphisme local pr{\`e}s de
tout point $ (x_1^\circ, x_{\frac{1}{2}}^\circ, x_0^\circ)$ tel que
$x_1^\circ$ soit inversible comme {\'e}l{\'e}ment de $J(c,1)$.
\end{proposition}
\proof Soit $\Phi$ la repr{\'e}sentation de l'alg{\`e}bre de
Jordan $J(c,1)$ dans $\text{End}(J(c,\frac{1}{2}))$ donn{\'e}e par (\ref{representation_J_c_1}).
D'apr{\`e}s les formules (\ref{decomp_frobenius1}) on a, en
posant, $\tau_c(x_{\frac{1}{2}})(x_1+x_0)=y_1+y_{\frac{1}{2}}+y_0$,
$$\begin{array}{l}
y_1=x_1,\\
y_{\frac{1}{2}}=\Phi(x_1)x_{\frac{1}{2}},\\
y_0=2L(e-c)L(x_{\frac{1}{2}})^2x_1+x_0.
\end{array}$$
Si l'{\'e}l{\'e}ment
$x_1^\circ$ consid{\'e}r{\'e} comme {\'e}l{\'e}ment de la sous-alg{\`e}bre
$J(c,1)$ est inversible,  l'endomorphisme $\Phi(x_1^\circ)$
est inversible. Ceci reste vrai pour $x_1$ voisin de $x_1^\circ$. Il en
r{\'e}sulte facilement que l'application consid{\'e}r{\'e}e est un
diff{\'e}omorphisme local pr{\`e}s de $ (x_1^\circ, x_{\frac{1}{2}}^\circ, x_0^\circ)$.
\hfill\hfill$\square$\\

Un {\'e}l{\'e}ment $x$ de $J$ est dit de {\it corang} $k$ si son rang est
$r-k$. Pour tout entier $k$ tel que $0\leq k\leq r$, on note dans la suite
$J_k$ l'ensemble des {\'e}l{\'e}ments de $J$ de corang $k$,
$$J_k=\{x\in J\tq \mu(x,0)=k\}.$$
Ainsi,  $J_0$ n'est
autre que l'ensemble $J^\times$ des {\'e}l{\'e}ments inversibles de $J$, et $J_r= \{0\}$.

\begin{proposition}\label{lemme_fondamental}
Soit   $k$ un entier tel que $1\leq k\leq r-1$. Alors $J_k$  est
une sous-vari{\'e}t{\'e} de $J$ de codimension $k+\frac{k(k-1)}{2}d$. En
particulier, Si $k=1$, alors $J_1$ est de codimension $1$, et si $ k\geq 2$, alors la codimension de $J_k$ est sup{\'e}rieur ou {\'e}gale {\`a} $3$. 
\end{proposition}
 \proof
Soit $1\leq k \leq r-1$. Soit $x^\circ\in J_k$ un {\'e}l{\'e}ment de rang $r-k$. Il
existe donc une d{\'e}composition de Peirce $e=\sum_{j=1}^r c_j$ telle
que 
$$ x^\circ=\sum_{j=1}^{r-k}\lambda_jc_j$$
avec $\lambda_j\neq 0$ pour tout $j$, $1\leq j\leq r-k$. Posons $c=\sum_{j=1}^{r-k}c_j$. L'{\'e}l{\'e}ment $c$ ainsi d{\'e}fini est un idempotent de
l'alg{\`e}bre de Jordan $J$. Soit 
$$J=J(c,1)\oplus J(c,\frac{1}{2})\oplus J(c,0)$$
la d{\'e}composition de Peirce de $J$ relative {\`a} cet idempotent. L'{\'e}l{\'e}ment
$x^\circ$ appartient {\`a} $J(c,1)$ et est un {\'e}l{\'e}ment inversible de $J(c,1)$.
Par suite, on peut appliquer la proposition \ref{lemme_Frobenius}. On obtient
ainsi des coordonn{\'e}es locales $(y_1(x),y_{\frac{1}{2}}(x), y_0(x))$ pr{\`e}s de
$(x^\circ,0,0)$, telles que $ x=\tau_c(y_{\frac{1}{2}}(x))(y_1(x)+y_0(x))$. L'action du groupe $L$ pr{\'e}servant la notion de rang d'un {\'e}l{\'e}ment, 
le rang de
$x$ est le m{\^e}me que le rang de
$y_1(x)+y_0(x)$. Mais le rang de $y_1+y_0$ est la somme du rang de
l'{\'e}l{\'e}ment $y_1$ et de l'{\'e}l{\'e}ment
$y_0$, pour tout $y_1\in J(c,1), y_0\in J(c,0)$. Pour $x$ voisin de
$x^\circ$, $y_1(x)$ est voisin de $x^{0}$ et donc inversible dans l'alg{\`e}bre de Jordan $J(c,1)$, ou si l'on
pr{\'e}f{\`e}re de rang $r-k$. Par cons{\'e}quent la condition 
$$ y_0(x) = 0$$
fournit un syst{\`e}me d'{\'e}quations locales de $J_k$ pr{\`e}s de
$x^\circ$. Ceci montre que $J_k$ est une sous-vari{\'e}t{\'e} de
codimension {\'e}gale {\`a} $\dim J(c,0)$. D'apr{\`e}s
(\ref{dim_J_0}) cette dimension est 
\begin{equation*}
\dim J(c,0)=k+\frac{k(k-1)}{2}d,
\end{equation*}
puisque l'idempotent $c$ est de rang $r-k$.\\ 
Si $k=1$ alors $\dim J(c,0)=1$ et  $J_1$ est de codimension
1. \\
Si $k\geq 2$, alors $\dim J(c,0)\geq 3$, puisque l'invariant de Peirce $d$
d'une alg{\`e}bre de Jordan euclidienne simple est non nul (voir \cite[page
70]{Faraut-Koranyi}), et donc la codimension de $J_k$ est sup{\'e}rieure ou {\'e}gale
{\`a} 3. \hfill\hfill$\square$\\

D'apr{\`e}s \cite[Proposition II.4.2]{Faraut-Koranyi}, le gradient de la
fonction $\det$ est polynomial et  v{\'e}rifie sur $J^\times$ la formule  
\begin{equation*}
\nabla\det(x)=\det(x)x^{-1}.
\end{equation*}
Plus explicitement, si $x=\sum_{j=1}^r\lambda_jc_j$ est une d{\'e}composition
de Peirce d'un {\'e}l{\'e}ment
inversible de $J$, alors 
\begin{equation}\label{grad_expl}
\nabla\det(x)=\sum_{j=1}^r(\lambda_1\cdots\overset{\lor}{\lambda}_j\cdots \lambda_r)c_j.
\end{equation}
L'ensemble $J\setminus J^\times=\sqcup_{k=1}^rJ_k$ est l'ensemble
des {\'e}l{\'e}ments $x\in J$ tels que $\det(x)=0$. C'est une vari{\'e}t{\'e} (au sens
alg{\'e}brique) de
codimension 1, dont on va maintenant pr{\'e}ciser le lieu singulier.

\begin{proposition}\label{lieu_Singulier_J} 
Soit $x\in\sqcup_{k=1}^rJ_k$. Alors, $\nabla\det(x)\not=0$ si, et seulement si, $x\in J_1$.
\end{proposition}
\proof Soit $x\in J_1$, alors il existe un rep{\`e}re de Jordan
$(c_j)_{1\leq j\leq r}$ et des nombres r{\'e}els non nuls $(\lambda_j)_{1\leq
  j\leq r-1}$ tels que $x=\sum_{j=1}^{r-1}\lambda_jc_j$. On pose
$x_\epsilon=\sum_{j=1}^{r-1}\lambda_jc_j+\epsilon c_r$ pour $\epsilon >0$
assez petit. D'apr{\`e}s la formule (\ref{grad_expl}) on a 
\begin{equation*}
\nabla
\det(x_\epsilon)=\sum_{j=1}^{r-1}\epsilon(\lambda_1\cdots\overset{\lor}{\lambda}_j\cdots
\lambda_{r-1}) c_j + (\lambda_1\cdots \lambda_{r-1})c_r.
\end{equation*}
Par cons{\'e}quent, lorsque $\epsilon \to 0$, $\nabla\det(x_\epsilon)\to \lambda_1\cdots
\lambda_{r-1}c_r$. Par continuit{\'e}, on a donc
\begin{equation*}
\nabla\det(x)=(\lambda_1\cdots \lambda_{r-1})c_r\not=0.
\end{equation*}
Si par contre $x\in \sqcup_{k=2}^r J_k$, un calcul similaire montre que
$\nabla\det(x)=0$. D'o{\`u} la proposition.

\hfill\hfill$\square$

\subsection{La stratification du cycle de Maslov}
{\'E}tant donn{\'e} un {\'e}l{\'e}ment $\sigma_0\in S$, et $k$ un entier, $0\leq k\leq r$,
on pose
\begin{equation*}
S_k(\sigma_0)=\{\sigma\in S \tq \mu(\sigma,\sigma_0)=k\}.
\end{equation*}
On a
\begin{equation*}
\Sigma(\sigma_0)=\bigsqcup_{1\leq k\leq r}S_k(\sigma_0).
\end{equation*}

Dans le cas o{\`u} $S$ est la vari{\'e}t{\'e} lagrangienne $\Lambda(r)$, pour un
sous-espace 
lagrangien $\lambda_0$ fix{\'e}, $S_k(\lambda_0)$ co{\"\i}ncide avec l'ensemble
$\Lambda^k(r)$ des sous-espaces lagrangiens $\lambda$ tels que
\text{dim}($\lambda\cap\lambda_0)=k$.

\begin{proposition}
$S_k(\sigma_0)$ est une sous-vari{\'e}t{\'e}  de $S$ de codimension $k+\frac{k(k-1)}{2}d$.
\end{proposition}
\proof
Soit $\sigma\in S_k(\sigma_0)$. Gr{\^a}ce {\`a} l'action du groupe $G$, on peut
toujours supposer que $\sigma_0=-e$, et $\sigma\top e$. Par la transformation
de Cayley, on est ramen{\'e} {\`a} l'{\'e}nonc{\'e} {\'e}quivalent pour $J_k$, qui a
{\'e}t{\'e} d{\'e}montr{\'e} pr{\'e}c{\'e}demment (proposition \ref{lemme_fondamental}).

\hfill\hfill$\square$

\begin{proposition} Le cycle de Maslov $\Sigma(\sigma_0)$ est
  une sous-vari{\'e}t{\'e} stratifi{\'e}e de $S$ (au sens alg{\'e}brique) de codimension 1.
  Le lieu singulier de $\Sigma(\sigma_0)$ est de codimension  sup{\'e}rieure ou {\'e}gale {\`a} $3$ dans $S$.
\end{proposition}
\proof 
Utilisant comme pr{\'e}c{\'e}demment la transformation de Cayley, on voir
d'abord, gr{\^a}ce {\`a} la proposition \ref{lieu_Singulier_J} que $S_1(\sigma_0)$ est
l'ensemble des points r{\'e}guliers de $\Sigma(\sigma_0)$, et $S_1(\sigma_0)$
est de codimension 1 dans $S$. Le lieu singulier est {\'e}gal {\`a}
$\sqcup_{2\leq k\leq r}S_k(\sigma_0)$, qui est de codimension $\geq 3$ dans
$S$, d'apr{\`e}s la proposition \ref{lemme_fondamental}.

\hfill\hfill$\square$

\subsection{L'orientation transverse du cycle de Maslov} 
Dans le cas de la vari{\'e}t{\'e} lagrangienne $\Lambda(r)$ (voir \cite{Arnold}),
Arnold montre que le cycle de
Maslov admet un champ continu de vecteurs
tangents {\`a} $\Lambda(r)$ et transverses {\`a}
$\Lambda^1(r)$. Un tel champs de vecteurs est
constitu{\'e}  par les vecteurs $v(\sigma)=\frac{d}{d\theta}(e^{i\theta}\sigma)_{|_{\theta=0}}$ pour
$\sigma\in\Lambda^1(r)$ et d{\'e}fini une orientation transverse. Nous allons montrer comment cette notion se
g{\'e}n{\'e}ralise {\`a} $S$ en utilisant sa structure causale.

\begin{proposition}
Soit $\sigma_0\in S$, et soit $\sigma\in S_1(\sigma_0)\subset
\Sigma(\sigma_0)$. Soit $H_{\sigma_0}(\sigma)$ l'espace tangent {\`a} $\Sigma(\sigma_0)$ au point
$\sigma$. Les vecteurs tangents aux courbes causales partant de $\sigma$ sont tous contenus dans un
m{\^e}me demi-espace ouvert de l'espace tangent {\`a} $S$ en $\sigma$ limit{\'e} par $H_{\sigma_0}(\sigma)$.
\end{proposition}
Par la transform{\'e}e de Cayley, on est ramen{\'e} {\`a} l'{\'e}nonc{\'e} {\'e}quivalent suivant.
\begin{proposition}
Soit $x\in J_1\subset J$. Soit $H(x)$ l'hyperplan tangent {\`a} $J_1$ en
$x$. Les vecteurs tangents aux courbes causales partant de $x$ sont tous contenus dans un m{\^e}me
demi-espace ouvert limit{\'e} par $H(x)$.
\end{proposition}
\proof
On peut supposer que $x= \sum_{j=1}^{r-1}\lambda_j c_j$, o{\`u}
$(c_j)_{1\leq j\leq r}$ est un rep{\`e}re de Jordan de $J$ et, pour $1\leq j\leq r-1$,  les $\lambda_j$
sont des scalaires non nuls. On sait alors qu'un vecteur normal (non nul !) est donn{\'e} par
$$ \nabla \det (x) = \big(\prod_{j=1}^{r-1}\lambda_j\big)\ c_r\ .$$
Si $y(t)$ est une courbe causale telle que $y(0)=x$ et $\dot{y} (0)=\omega\in \Omega$, on a donc
$$ \frac{d}{dt}(\det y(t))_{\mid_{t=0}} = \langle \dot{y}(0)|\nabla\det x\rangle 
 = \bigl(\prod_{j=1}^{r-1}\lambda_j\bigr)\, \langle \omega|c_r \rangle.$$
Comme $\omega\in \Omega$ et $c_r\in \overline{\Omega} \setminus \{0\}$, on a $\langle\omega|c_r\rangle>0 $, et
par suite le signe de $\langle\dot{y}(0)|\nabla\det x\rangle$ ne d{\'e}pend pas
de la courbe causale choisie. 
\hfil\hfill$\square$\\

{\`A} tout point $\sigma\in S_1(\sigma_0)$, on associe le demi-espace ouvert
$H^+_{\sigma_0}(\sigma)$ contenant les vecteurs tangents aux courbes
causales passant par $\sigma$. La famille de ces sous-espaces d{\'e}finit une
{\it orientation transverse\/} canonique du cycle de Maslov $\Sigma(\sigma_0)$. De plus ces orientations transverses sont compatibles {\`a} l'action de $G$, en ce
sens que si $g\in G$, alors l'orientation transverse canonique du cycle de
Maslov $\Sigma(g\sigma_0)$ au point $g\sigma\in S_1(g\sigma_0)$ est donn{\'e}e
par le demi-espace $ Dg(\sigma) \big(H^+_{\sigma_0}(\sigma)\big)=
H^+_{g\sigma_0}(g\sigma)$, puisque la notion de courbe causale est stable par l'action de $G$.

\section{Le nombre d'Arnold, l'indice d'inertie et l'indice  d'Arnold-Leray-Maslov }
Les principales propri{\'e}t{\'e}s topologiques des cycles de Maslov d{\'e}montr{\'e}es
dans le paragraphe pr{\'e}c{\'e}dent g{\'e}n{\'e}ralisent celles montr{\'e}es par Arnold
dans le cas de la  vari{\'e}t{\'e} lagrangienne $\Lambda(r)$. Dans ce cas particulier,
Arnold (voir \cite{Arnold} et \cite{Arnold2}) utilise ces propri{\'e}t{\'e}s pour
construire un indice pour les couples de  points transverses du rev{\^e}tement
universel $\widetilde{\Lambda(r)}$. Les arguments d'Arnold sont encore valables dans
notre situation et nous permettent dans la suite de proposer une deuxi{\`e}me
construction, g{\'e}om{\'e}trique cette fois-ci, d'une primitive.

Soit $\sigma_0\in S$, et soit $\gamma$ un chemin dans $S$, d'origine
$\gamma(0)=\sigma$, et d'extr{\'e}mit{\'e}
$\gamma(1)=\tau$. Un tel chemin est dit {\it admissible\/} (relativement au cycle de Maslov
$\Sigma_0=\Sigma(\sigma_0)$) s'il est de classe $\mathcal{C}^1$ par morceaux, si ses extr{\'e}mit{\'e}s
$\sigma$ et $\tau$ n'appartiennent pas {\`a}
$\Sigma_0$, si le chemin n'intersecte le cycle de Maslov que pour un nombre fini de valeurs de
$t\in [0,1]$, (disons $ t=t_1,t_2,\dots, t_k$), si chaque point d'intersection $\gamma(t_j)$
est un point r{\'e}gulier de $\Sigma_0$ (c'est-{\`a}-dire $\gamma(t_j)\in S_1(\sigma_0))$, et enfin si en un tel point, la courbe $\gamma$ est
diff{\'e}rentiable et transverse {\`a}
$\Sigma_0$. {\it Le nombre d'Arnold\/} d'un tel chemin admissible est par
d{\'e}finition la somme
$\epsilon_1+\epsilon_2+\dots+\epsilon_k$, o{\`u} $\epsilon_j$ vaut $+1$ si $\dot \gamma (t_j)\in
H^+_{\sigma_0}(\gamma(t_j))$, et $\epsilon_j$ vaut $-1$ dans le cas
contraire.

\begin{theoreme}
Soient $\sigma$ et $\tau$ deux points n'appartenant pas au cycle de
Maslov $\Sigma_0$. Toute classe d'homotopie de chemins d'origine $\sigma$ et
d'extr{\'e}mit{\'e} $\tau$ contient un chemin admissible. Deux chemins admissibles
d'origine $\sigma$ et d'extr{\'e}mit{\'e} $\tau$ qui sont homotopes ont m{\^e}me
nombre d'Arnold.
\end{theoreme}
La d{\'e}monstration est identique au cas classique (voir
\cite{Arnold2} ou \cite{Cappell-Lee-Miller}), au vu des r{\'e}sultats du
paragraphe 7. \\
Ce th{\'e}or{\`e}me permet d'{\'e}tendre la notion de nombre d'Arnold {\`a}
un chemin quelconque $\gamma$ d'origine $\sigma$ et d'extr{\'e}mit{\'e} $\tau$, en le
d{\'e}clarant {\'e}gal par d{\'e}finition {\`a} l'indice d'un chemin admissible de m{\^e}mes
origine et extr{\'e}mit{\'e} et homotope {\`a}
$\gamma$.\\

 On va maintenant tirer profit de l'invariance par homotopie du
nombre d'Arnold pour d{\'e}finir l'{\it indice d'Arnold\/} d'un couple de
points du rev{\^e}tement universel 
$\widetilde S$.\\
Soient $\widetilde \sigma_0$ et $\widetilde \tau_0$ deux points de $\widetilde S$, de projections
respectives $\sigma_0$ et $\tau_0$ sur $S$. Soit
$\sigma(t)$, $0\leq t\leq 1$ une courbe causale d'origine $\sigma_0$, et $ \tau
(t)$, $0\leq t\leq 1$ une courbe anti-causale d'origine $\tau_0$. On note $ \widetilde{\sigma(t)}$ (resp. $\widetilde {\tau(t)} $) la courbe dans $\widetilde S$ d'origine $\widetilde{\sigma_0}$  (resp. $ \widetilde {\tau_0} $) relevant
$\sigma(t)$ (resp. $ \tau(t)$). D'apr{\`e}s la proposition
\ref{estim_det_c_causale_S}, il existe $\epsilon >0$ tel que pour tout $t,
0<t<\epsilon$, les points $\sigma(t)$ et $\tau(t)$ n'appartiennent pas au cycle de Maslov $\Sigma_0=\Sigma(\sigma_0)$.
Fixons un tel $t$, et soit  $\gamma_t(s)$, $0\leq s\leq 1$ un chemin admissible
(relativement {\`a} $\Sigma_0$) d'origine $\sigma(t)$ et d'extr{\'e}mit{\'e}
$\tau(t)$, dont le relev{\'e} dans
$\widetilde S$ a pour origine $\widetilde{\sigma}(t)$ et pour extr{\'e}mit{\'e}
$\widetilde{\tau}(t)$.\\

\vskip4.7cm

\begin{pspicture}(0,0)(0,0)
\pscurve[linecolor=blue,showpoints=true](2,2.5)(4,3.5)(5,2.75)(7.5,4.5)(9,3)(9.7,3.5)
\pscurve[linecolor=green, showpoints=true](1,1)(1.5,2)(2,2.5)
\pscurve[linecolor=green, showpoints=true](9.7,3.5)(10,4)(10.5,5)
\psellipse[linecolor=red](1,2)(1,0.25)
\psline[linecolor=red](2,2)(1,1)(0,2)
\psellipse[linecolor=red](10.5,4)(1,0.25)
\psline[linecolor=red](9.5,4)(10.5,5)(11.5,4)
\rput(10.5,5.3){$\tau_0$}
\rput(1,0.77){$\sigma_0$}
\rput(1,2){$\sigma(t')$}
\rput(2.5,2.5){$\sigma(t)$}
\rput(4,3.8){$\gamma_t(s_1)$}
\rput(5,2.5){$\gamma_t(s_2)$}
\rput(7.5,4.75){$\gamma_t(s_3)$}
\rput(9,2.8){$\gamma_t(s_k)$}
\rput(9.3,3.5){$\tau(t)$}
\rput(10.5,4){$\tau(t')$}
\end{pspicture}
\begin{center}
Fig. 1
\end{center}
\medskip
Par d{\'e}finition {\it l'indice d'Arnold\/}
$\nu(\widetilde{\sigma}_0,\widetilde{\tau}_0)$ est le nombre d'Arnold du chemin
$\gamma_t$. Clairement, {\`a} $t$ fix{\'e}, cet indice ne d{\'e}pend pas du choix du
chemin $\gamma_t$, puisque tout autre chemin ayant les m{\^e}mes propri{\'e}t{\'e}s lui est homotope et donc a m{\^e}me nombre
d'Arnold.

De plus,  l'indice est ind{\'e}pendant du param{\`e}tre $t$ (voir Fig 1).  En effet, soit
$t'$ une autre valeur, avec toujours
$0<t'<\epsilon$. On peut toujours supposer que $t'<t$. On consid{\`e}re le
chemin d'origine $ \sigma(t')$ et d'extr{\'e}mit{\'e} $\tau(t')$ obtenu en concat{\'e}nant
le chemin $\big(\sigma(u)$, $t'\leq u\leq t\big)$, le chemin $\gamma_t$ de
$\sigma(t)$ {\`a} $\tau(t)$, et enfin le chemin $\big(\tau(v)$, $t'\leq v \leq t)$.
Il est imm{\'e}diat que le nombre d'Arnold de ce chemin concat{\'e}n{\'e} est le
m{\^e}me que le nombre d'Arnold du chemin $\gamma_t$, puisque les deux
bouts de chemin concat{\'e}n{\'e}s {\`a} $\gamma_t$ ne rencontrent pas le cycle de
Maslov $\Sigma_0$, comme on le d{\'e}duit ais{\'e}ment de  la proposition
\ref{estim_det_c_causle}.

Enfin, cet indice est ind{\'e}pendant des courbes (causale ou
anti-causale) $\sigma(t)$ et $\tau(t)$. En effet, {\'e}tant donn{\'e} deux courbes causales
$ \sigma_1(t),0\leq t\leq 1$ et $\sigma_2(t), 0\leq t\leq 1$ de m{\^e}me origine
$\sigma_0$, on peut pour $t$ assez petit, joindre $\sigma_1(t)$ {\`a} $\sigma_2(t)$ par
une courbe qui ne rencontre pas le cycle de Maslov $\Sigma_0$, comme on le voit en
utilisant la convexit{\'e} du c{\^o}ne $\Omega$.
\vskip4.9cm
\hskip3cm
\rotatebox{-45}{
\begin{pspicture}(1,1)(1,1)
\pscurve[linecolor=green](1,1)(0.5,2.5)(0,4)(-0.5,6.5)
\pscurve[linecolor=green](1,1)(1.5,2.5)(2,4)(2.5,6.5)
\pscurve[linecolor=gray,linestyle=dashed]{->}(0,4)(1,4.1)(2,4)
\psellipse[linecolor=red](1,5.5)(2,0.5)
\psline[linecolor=red](-1,5.5)(1,1)(3,5.5)
\end{pspicture}
}
\begin{pspicture}(1,1)(1,1)
\rput(0.8,0.7){$\sigma_0$}
\rput(4,5.6){$\sigma_1(t)$} 
\rput(6,3.5){$\sigma_2(t)$} 
\end{pspicture}
\begin{center}
Fig. 2
\end{center}

\bigskip

Un argument analogue vaut pour la courbe anti-causale d'origine $\tau_0$. Notons
que si $\tau_0$ est transverse {\`a} $\sigma_0$, il n'est pas n{\'e}cessaire d'utiliser
une telle courbe anti-causale, et on peut supposer que le chemin $\gamma_t$ est
d'extr{\'e}mit{\'e} $\tau_0$.\\

La construction de l'indice d'Arnold ne fait appel qu'{\`a}
des notions invariantes par transformations causales (courbes causales, cycle de
Maslov, transversalit{\'e}), et donc l'indice d'Arnold est invariant par l'action de
$\Gamma$.\\

On va maintenant calculer l'indice d'Arnold ``en coordonn{\'e}es'', comme nous
l'avons fait pour l'indice de Souriau. On fixe pour cela un rep{\`e}re de
Jordan $(c_j)_{ 1\leq j\leq r}$ de $J$. Gr{\^a}ce aux r{\'e}sultats concernant les
orbites de l'action de $G$ dans $S\times S$ (voir paragraphe 4), la
proposition suivante couvre le cas g{\'e}n{\'e}ral.

\begin{proposition}\label{ind_Maslov-Arnold_coor}
 Soient les deux points de $\widetilde{S}$, 
$$\widetilde{\sigma}_0=\widetilde{-e} =(-e,-\pi)\;\;\text{et}\;\;
\widetilde{\tau}_0=\big(-\sum_{j=1}^\ell c_j+\sum_{j=\ell+1}^r
e^{i\varphi_j}c_j,\ \varphi\big).$$
On note $\sigma_0=-e$ et $\tau_0$ leurs projections correspondantes dans $S$. On
suppose que
\begin{enumerate}
\item[$(i)$] $-\pi<\varphi_j<\pi$, $\forall j$, $\ell+1\leq j\leq r$;
\item[$(ii)$] $r\varphi= -\ell\pi+\sum_{j=\ell+1}^r\varphi_j+2k\pi$,
  avec $k\in \mathbb{Z}.$
\end{enumerate}
Alors
\begin{equation}\label{Indice_Arnold_en_Coor}
\nu(\widetilde{\sigma}_0,\widetilde{\tau}_0)
=-\ell+k=k-\mu(\sigma_0,\tau_0).
\end{equation}
\end{proposition}
\proof Le point
$\sigma=\sum_{j=1}^re^{i\theta_j}c_j$ appartient au
cycle de Maslov $\Sigma_0=\Sigma(-e)$ si et seulement s'il existe
$j$, $1\leq j\leq r$ tel que $\theta_j\in -\pi+2\pi\mathbb{Z}$. Il est
r{\'e}gulier (c'est-{\`a}-dire $\sigma\in S_1(-e)$) si et seulement si $
\theta_j\in -\pi+2\pi\mathbb{Z}$ pour une seule valeur de l'indice $j$. Si, par exemple, $\sigma=-c_1+\sum_{j=2}^re^{i\theta_j}c_j$,
avec pour $2\leq j\leq r$, $\theta_j\notin -\pi+2\pi\mathbb{Z}$,  alors un
vecteur normal au plan tangent en $\sigma$ {\`a} $\Sigma_0$ donnant l'orientation positive est le vecteur $-ic_1$.\\

Comme courbe causale d'origine $\sigma_0$, choisissons la courbe
$$ \sigma(t)=\sum_{j=1}^re^{i(-\pi+t)}c_j$$
dont le relev{\'e} dans $\widetilde{S}$ partant de $\widetilde{\sigma}_0$ est
$$\widetilde{\sigma(t)} =\big(\sum_{j=1}^re^{i(-\pi+t)}c_j, -\pi+t)\ ,$$
et pour courbe anti-causale d'origine $\widetilde{\tau}_0$
$$\tau(t) = \sum_{j=1}^\ell e^{i(-\pi-t)}c_j+\sum_{j=\ell+1}^re^{i\varphi_j}c_j$$
dont le relev{\'e} dans $\widetilde{S}$ est
$$ \widetilde{\tau(t)}=\Big(\tau(t),
\frac{1}{r}\big[\ell(-\pi-t)+\sum_{j=\ell+1}^r\varphi_j+2k\pi\big]\Big).$$

On notera que la courbe $\tau(t)$ n'est pas, au sens strict, anti-causale,
puisque son vecteur tangent en chaque point est seulement dans la fronti{\`e}re du
c{\^o}ne causal. N{\'e}anmoins, il est vrai que pour
$t>0$ assez petit, le point $ \tau(t)$ n'est pas dans $\Sigma_0$, et permet
donc le calcul de l'indice d'Arnold du couple $
(\widetilde{\sigma}_0,\widetilde{\tau}_0)$.\\

Ayant fix{\'e} un tel $t>0$ assez petit, le chemin  $\gamma_t$ va {\^e}tre obtenu par
concat{\'e}nation de plusieurs morceaux de chemins admissibles, dont on va calculer
pour chacun la contribution {\`a} l'indice de d'Arnold. Pour la commodit{\'e} de la lecture, chaque morceau de
chemin  est param{\'e}tr{\'e} par $s\in[0,1]$.\\

Le premier morceau de chemin est donn{\'e} par
$$ \gamma_t(s)= e^{i[(1-s)(-\pi+t)+s(-\pi-t)]}c_1+ \sum_{j=2}^re^{is(-\pi+t)}
c_j.$$
Pour $s=\frac{1}{2}$, $\gamma_t(s)=-c_1 +\sum_{j=2}^re^{it(-\pi+t)}
c_j$, et c'est le seul point d'intersection du chemin $\gamma_t$ avec
$\Sigma_0$. De plus,$\dot{\gamma}_t(\frac{1}{2})=2itc_1$, de sorte que  le
chemin $\gamma_t$ coupe donc transversalement le cycle de Maslov, et la contribution {\`a}
l'indice est $-1$.\\

On continue par le m{\^e}me type de  chemin, pour chaque $j$, $2\leq j\leq \ell$.
L'extr{\'e}mit{\'e} finale de ces chemins est alors le point
$$\sigma'_t= \sum_{j=1}^\ell e^{i(-\pi-t)} c_j+\sum_{j=\ell+1}^r e^{i(-\pi+t)}c_j$$ et
la contribution {\`a} l'indice d'Arnold est $-\ell$.\\

On consid{\`e}re maintenant le chemin
$$ \gamma_t(s) = \sum_{j=1}^\ell e^{i(-\pi-t)}
c_j+\sum_{j=\ell+1}^re^{((1-s)(-\pi+t) + s\varphi_j)}c_j.$$
Ce chemin ne rencontre pas le cycle de Maslov $\Sigma_0$, et donc la
contribution de ce morceau de chemin {\`a} l'indice est  nulle. Le point d'arriv{\'e}e
est maintenant
$$ \tau'_t=\sum_{j=1}^\ell e^{i(-\pi-t)}c_j+ \sum_{j=\ell+1}^r e^{i\varphi_j}c_j\ .$$

Un calcul facile permet de d{\'e}terminer l'extr{\'e}mit{\'e} du relev{\'e} dans
$\widetilde{S}$ du chemin parcouru jusqu'{\`a} pr{\'e}sent. On trouve que c'est le point
$$\widetilde{\tau'_t} = \big(\tau'_t, \frac{1}{r}[\ell(-\pi-t)+\sum_{j=\ell+1}^r\varphi_j]\big).$$
Pour terminer, on a besoin du lemme suivant.
 
\begin{lemme}
Soient donn{\'e}s $(\psi_j)_{1\leq j\leq r}$ $r$ nombres
r{\'e}els tels que $\psi_j\notin -\pi+2\pi\mathbb{Z}$. Soit
$\gamma
$ le chemin d{\'e}fini pour
$0\leq s\leq 1$ par
$$ \gamma(s) = e^{i(\psi_1+2s\pi)}c_1+\sum_{j=2}^r e^{i\psi_j}c_j$$
Le chemin $\gamma$ est admissible relativement au cycle de Maslov $\Sigma_0$, et
son nombre d'Arnold vaut $+1$.
\end{lemme}
En effet, le chemin $\gamma$ ne coupe $\Sigma_0$ que pour une seule valeur de
$s$, transversalement, et dans le sens positif. Donc le nombre
d'Arnold de ce chemin vaut
$+1$.

Ce lemme permet de calculer l'indice d'Arnold de deux points de $\widetilde S$
ayant m{\^e}me projection dans $S$.  En particulier, en passant du
point
$\widetilde {\tau'_t}$ au point
$\widetilde {\tau_t}$, la contribution {\`a} l'indice d'Arnold est
$k$. D'o{\`u} la proposition \ref{ind_Maslov-Arnold_coor}.

\hfill\hfill$\square$

\begin{corollaire}\label{pro_lien_m_nu} Soient
  $\widetilde{\sigma},\widetilde{\tau}$ deux {\'e}l{\'e}ments
  de  $\widetilde{S}$ de projections respectives $\sigma$ et $\tau$. Alors
\begin{equation}\label{lien_m_nu}
\nu(\widetilde{\sigma},\widetilde{\tau})=\frac{1}{2}\bigl(m(\widetilde{\sigma},\widetilde{\tau})-\mu(\sigma,\tau)-r\bigr).
\end{equation}
\end{corollaire}
\proof 
Il suffit en effet de comparer les formules (\ref{eq:m_en_coord_simple}) et
(\ref{Indice_Arnold_en_Coor}) et d'utiliser l'invariance des deux membres
par $\Gamma$.

\hfill\hfill$\square$

\begin{remarque}{\rm
Un sous-produit de ce corollaire est que le membre de droite de la formule (\ref{lien_m_nu}) est un entier (et non un demi-entier). On tire profit de
ce r{\'e}sultat  en introduisant la notion d'{\it indice d'inertie\/} et d'{\it indice d'Arnold-Leray-Maslov\/}.
}\end{remarque}
 
\begin{definition}
  \begin{itemize}
  \item[$(i)$] Soient $(\sigma_1,\sigma_2,\sigma_3)\in S^3$. On
appelle indice d'inertie du triplet $(\sigma_1,\sigma_2,\sigma_3)$ la quantit{\'e}
\begin{equation*} 
\jmath(\sigma_1,\sigma_2,\sigma_3) = \frac{1}{2}\big(\iota(\sigma_1,\sigma_2,\sigma_3)+\mu(\sigma_1,\sigma_2) -\mu(\sigma_1,\sigma_3)+\mu(\sigma_2,\sigma_3)+r\big).
\end{equation*}
  \item[$(ii)$] Soient $\widetilde{\sigma}_1,\widetilde{\sigma}_2 \in
    \widetilde{S}$, de projections respectives $\sigma_1,\sigma_2$. On
    appelle indice d'Arnold-Leray-Maslov la quantit{\'e}
\begin{equation*}
n(\widetilde{\sigma}_1,\widetilde{\sigma}_2) = \nu(\widetilde{\sigma}_1,\widetilde{\sigma}_2)+\mu(\sigma_1,\sigma_2)+r.
\end{equation*}
  \end{itemize}
\end{definition}

\begin{theoreme} 
L'indice d'inertie satisfait les propri{\'e}t{\'e}s suivantes :
\begin{itemize}
\item[$(i)$] $\jmath$ est  {\`a} valeurs enti{\`e}res.
\item[$(ii)$] $\jmath$ v{\'e}rifie la relation de 2-cocycle \footnote{On
    observera que l'indice d'inertie $\jmath$ n'a pas la propri{\'e}t{\'e} d'antisym{\'e}trie
    que poss{\`e}de  l'indice de Maslov $\imath$.}
\begin{equation*}
\jmath(\sigma_1,\sigma_2,\sigma_3)-\jmath(\sigma_1,\sigma_2,\sigma_4)+\jmath(\sigma_1,\sigma_3,\sigma_4)-\jmath(\sigma_2,\sigma_3,\sigma_4)=0
\end{equation*}
pour tous $\sigma_1,\sigma_2,\sigma_3\in S$.
\item[$(iii)$] l'indice d'Arnold-Leray-Maslov est une primitive de
l'indice d'inertie, c'est-{\`a}-dire
\begin{equation*}
\jmath(\sigma_1,\sigma_2,\sigma_3) =n(\widetilde{\sigma}_1,\widetilde{\sigma}_2)-n(\widetilde{\sigma}_1,\widetilde{\sigma}_3)+n(\widetilde{\sigma}_2,\widetilde{\sigma}_3)
\end{equation*}
pour tous $\widetilde{\sigma}_1,\widetilde{\sigma}_2,\widetilde{\sigma}_3\in \widetilde{S}$, de
projections respectives $\sigma_1,\sigma_2,\sigma_3$.
\end{itemize}
\end{theoreme}
\proof $(i)$ et $(ii)$ sont des cons{\'e}quences faciles de
$(iii)$. Pour $(iii)$, on utilise la relation
\begin{equation*}
n(\widetilde{\sigma},\widetilde{\tau}) = \frac{1}{2}\bigl(m(\widetilde{\sigma},\widetilde{\tau})+\mu(\sigma,\tau)+r\bigr)
\end{equation*}
valable pour tous $\widetilde{\sigma},\widetilde{\tau} \in \widetilde{S}$ de
projections respectives $\sigma,\tau$, qui est cons{\'e}quence du
corollaire \ref{pro_lien_m_nu}. On voit alors que la relation
$(iii)$ est {\'e}quivalente {\`a} la propri{\'e}t{\'e} de cocycle pour l'indice de
Maslov (\ref{rel_cocy_g}).

\hfill\hfill$\square$

\section{L'indice de Maslov d'une paire de chemins dans $S\times S$}
Dans ce paragraphe, nous allons adapter  une m{\'e}thode due {\`a} Cappell-Lee-Miller
\cite{Cappell-Lee-Miller} pour expliquer (sans plus de d{\'e}tails) comment on
peut d{\'e}finir l'indice de Maslov d'une paire de
chemins de $S$. \\

Un chemin $[0,1]\ni t\mapsto \varpi(t)=(\sigma_1(t),\sigma_2(t))\in S\times
S$ est {\it propre}, si $\sigma_1(t)$ est transverse {\`a}
$\sigma_2(t)$ pour $t=0$ et $t=1$.\\

 Soit $\varpi$ un tel chemin et supposons
qu'il est r{\'e}gulier. D'apr{\`e}s le paragraphe 7, si on note
\begin{equation*}
\mathcal{Z}=\{(t,\sigma)\in[-1,1]\times S , \; \det(\sigma_1(t)-\sigma)=0\}
\end{equation*}
 alors $\mathcal{Z}\cap \{t\}\times S=\{t\}\times \Sigma(\sigma_1(t))$ est
 sous-vari{\'e}t{\'e} stratifi{\'e}e (au sens alg{\'e}brique), transversalement orient{\'e}e, de
 codimension 1 de $\{t\}\times S$ et dont la singularit{\'e} est de
 codimension $\geq 3$ dans $\{t\}\times S$. Cette singularit{\'e} est donc de
 codimension $\geq 2$ dans $\mathcal{Z}$. Par cons{\'e}quent moyennant une petite
 perturbation du
 chemin $\varpi$, tout en gardant ses extr{\'e}mit{\'e}s fixes, on peut transformer
 le chemin $t\mapsto (t,\sigma_2(t))$ en un autre chemin $t\mapsto (t,\sigma_2'(t))$ ne
 rencontrant $\mathcal{Z}$ qu'en un nombre fini de points situ{\'e}s sur la
 strate sup{\'e}rieure $\{t\}\times S_1(\sigma_1(t))$ et ceci d'une mani{\`e}re transverse.\\
\begin{definition}
L'indice de Maslov $\nu_p(\varpi)$ du chemin propre $\varpi$
est par d{\'e}finition le nombre d'intersections, compt{\'e}s avec leurs signes, du
chemin $t\mapsto (t,\sigma_2'(t))$ avec $\mathcal{Z}$.
\end{definition}

 La r{\'e}union disjointe de toutes les  strates sup{\'e}rieures
\begin{equation*}
\bigcup_{0<t<1}\{t\}\times S_1(\sigma_1(t))
\end{equation*}
est une sous-vari{\'e}t{\'e} ouverte de $\mathcal{Z}$ dont la singularit{\'e}
est de codimension sup{\'e}rieure ou {\'e}gale {\`a} 3 dans $\mathcal{Z}$. On en
d{\'e}duit que deux petites perturbations du chemin $(\cdot,\sigma_2'(\cdot))$
sont homotopes et par cons{\'e}quent 
$\nu_p(\varpi)$ est ind{\'e}pendant du choix de cette perturbation.\\

Remarquons enfin que si le chemin $\varpi$ est seulement continue, on peut
l'approcher par un chemin r{\'e}gulier et d{\'e}finir  son indice
$\nu_p(\varpi)$. Comme la singularit{\'e} est de codimension $\geq 3$, l'indice
$\nu_p(\varpi)$ est ind{\'e}pendant du choix de cette approximation.\\

Les deux propositions suivantes 
 sont une cons{\'e}quence directe de la d{\'e}finition de l'indice
 $\nu_p$.

\begin{proposition}\label{pptes_mu_p} Soit $\varpi : [0,1]\to S\times S$, $t\mapsto
  (\sigma_1(t),\sigma_2(t))$ un chemin propre et 
  r{\'e}gulier.
  \begin{enumerate}
  \item Si $\sigma_1(t)$ est transverse {\`a} $\sigma_2(t)$ pour tout
    $t\in[0,1]$, alors $\nu_p(\varpi)=0$.
  \item Si $[0,1]\ni t\mapsto \hat{\varpi}(t)=(\sigma_1(-t),\sigma_2(-t))$
    est le chemin inverse de chemin $\varpi$, alors
    $\nu_p(\varpi)=\nu_p(\hat{\varpi})$.
  \item Pour tout $\alpha\in ]0,1[$, 
\begin{equation*}
\nu_p(\varpi)=\nu_p(\varpi_{|_{[0,\alpha]}})+\nu_p(\varpi_{|_{[\alpha,1]}})
\end{equation*}
  \item Soit $g=(g_t)_{0\leq t\leq 1}$ une famille d'{\'e}l{\'e}ments de $G$ telle que
$[0,1]\ni t\mapsto g_t$ soit continue. Si on pose $(g*\varpi)(t)=(g_t(\sigma_1(t)),g_t(\sigma_2(t))$ pour $t\in[0,1]$, alors
\begin{equation*}
\nu_p(g*\varpi)=\nu_p(\varpi)
\end{equation*}
  \end{enumerate}
\end{proposition}
\hfill\hfill$\square$

\begin{proposition}\label{chemins_homotopes}
 Si $[a,b]\times[0,1]\ni (s,t)\mapsto
    (\sigma_1(s,t),\sigma_2(s,t))\in S\times S$ est une application continue, et si pour
    tout $s\in[a,b]$, le chemin $\varpi_s(t)=(\sigma_1(s,t),\sigma_2(s,t))$
    est propre et r{\'e}gulier dans $S\times S$, alors
\begin{equation*}
\nu_p(\varpi_0)=\nu_p(\varpi_1).
\end{equation*}
\end{proposition}
\hfill\hfill$\square$

La fonction $f: \theta\mapsto \det(\sigma_1(t)-e^{i\theta}\sigma_2(t))$ a un
nombre fini de z{\'e}ros (modulo $2\pi$). Il existe donc $\epsilon\in]0,\pi[$
tels que les  $\theta$ qui v{\'e}rifient  $0<|\theta|<\pi$ ne soient pas des
z{\'e}ros de $f$. D'o{\`u} le lemme suivant : 
 
\begin{lemme}\label{lemme_de_perturbation}
Soit $\varpi : [0,1]\to S\times S$, $t\mapsto
  (\sigma_1(t),\sigma_2(t))$ un chemin r{\'e}gulier dans $S\times S$. Alors il
  existe $\epsilon\in ]0,\pi[$ tel que pour tout $\theta$ v{\'e}rifiant $0<|\theta|<\epsilon$, le nouveau chemin perturb{\'e}
  \begin{equation*}
\varpi_\theta(t)=(\sigma_1(t),e^{i\theta}\sigma_2(t))
\end{equation*}
 soit propre.
\end{lemme}
\hfill\hfill$\square$

La proposition suivante se d{\'e}duit facilement de  la proposition \ref{chemins_homotopes} et du
lemme \ref{lemme_de_perturbation}.

\begin{proposition}\label{mu_bien_defini} Soit $\varpi$ un chemin r{\'e}gulier dans $S\times S$. Si
  $\varpi_{\theta_1}$ et $\varpi_{\theta_2}$ sont deux perturbations
  diff{\'e}rentes de $\varpi$ en chemins propres, alors
\begin{equation*}
\nu_p(\varpi_{\theta_1})=\nu_p(\varpi_{\theta_2})
\end{equation*}
\end{proposition}
\hfill\hfill$\square$

\begin{definition}
Soit $\varpi$ un chemin r{\'e}gulier dans $S\times S$. L'indice de Maslov de
$\varpi$ est par d{\'e}finition l'indice de Maslov d'une perturbation
$\varpi_\theta$ de $\varpi$ :
\begin{equation*}
\nu(\varpi)=\nu_p(\varpi_\theta).
\end{equation*}
\end{definition}

La proposition \ref{mu_bien_defini} montre que l'indice
$\nu(\varpi)$ est bien d{\'e}fini et ne d{\'e}pend pas de la perturbation
choisie. 

Finalement les propositions \ref{pptes_mu_p} et \ref{chemins_homotopes} s'{\'e}tendent {\`a}
l'indice de Maslov $\nu(\varpi)$ d'un chemin quelconque $\varpi$. Pour plus
de d{\'e}tails nous invitons le lecteur {\`a} consulter le travail de
Cappell-Lee-Miller \cite{Cappell-Lee-Miller} dans lequel cette th{\'e}orie
est tr{\`e}s bien d{\'e}velopp{\'e}e.

\section{Le nombre de rotation g{\'e}n{\'e}ralis{\'e}}
Dans \cite{Barge-Ghys}, les auteurs montrent qu'on peut utiliser une primitive
(quelconque) du cocycle classique de Maslov pour construire un invariant
(pour la conjugaison) sur le groupe symplectique, baptis{\'e} {\it nombre de
rotation symplectique}. Nous montrons qu'une construction analogue est
possible dans la situation consid{\'e}r{\'e}e dans ce travail.

Commen{\c c}ons par quelques lemmes {\`a} caract{\`e}re g{\'e}n{\'e}ral.
\begin{lemme}
Soit $G$ un groupe de Lie semi-simple {\`a} centre fini
et connexe. Soit $G=KAN$ une d{\'e}composition d'Iwasawa de $G$, et soit
$r=\dim A$ le rang d{\'e}ploy{\'e} de $G$. Alors
\begin{itemize}
\item[$(i)$] tout {\'e}l{\'e}ment de $N$ est un commutateur,
\item[$(ii)$] tout {\'e}l{\'e}ment de $A$ peut s'{\'e}crire comme produit d'au plus $r$
commutateurs. 
\end{itemize}
\end{lemme}
\proof
$(i)$
Soit $\mathfrak{a} = Lie(A)$, $H$ un {\'e}l{\'e}ment
r{\'e}gulier de $\mathfrak{a}$ et $h=\exp H$.  Alors l'application
$$ n\mapsto hnh^{-1}n^{-1}, \quad N\to N$$ est un
diff{\'e}omorphisme de $N$ sur lui-m{\^e}me (voir \cite[Ch. VI. Corollary 4.8]{Helgason}).\\

$(ii)$ Soit $\Sigma$ le
syst{\`e}me de racines restreintes relativement
{\`a} $\mathfrak{a}$ et $\Sigma^+$ l'ensemble des
racines positives tel que $$ Lie(N) = \mathfrak{n}
= \bigoplus_{\alpha \in \Sigma^+} \mathfrak{g}_\alpha\ .$$
Soit $\{\alpha_1,\alpha_2,\dots,
\alpha_r\}$ les racines simples correspondantes. Pour $j, 1\leq
j\leq r$, on note $H_j$ le covecteur associ{\'e} {\`a} la racine $\alpha_j$. La
sym{\'e}trie orthogonale associ{\'e}e {\`a} la racine $\alpha_j$ est un {\'e}l{\'e}ment du
groupe de Weyl, et on en choisit un repr{\'e}sentant $m_j$ dans le
normalisateur de $\mathfrak{a}$ dans $K$. On a alors, pour tout $j, 1\leq j\leq
r$
\begin{equation}\label{commutateur}
m_j\exp H_jm_j^{-1}= \exp(- H_j).
\end{equation}

Maintenant, soit $h=\exp H$ un {\'e}l{\'e}ment quelconque de $A$. On a
$$ H=\sum_{j=1}^rt_jH_j,\quad h= \prod_{j=1}^r\exp t_jH_j,\quad  t_j\in \mathbb{R}, 1\leq j\leq
r\ .$$
D'apr{\`e}s (\ref{commutateur})
$$ \exp t_jH_j = m_j\exp(-\frac{t_j}{2}H_j)m_j^{-1}\exp(\frac{t_j}{2}H_j)$$
est donc un commutateur, et l'assertion  $(ii)$ en d{\'e}coule.

\hfill\hfill$\square$
\begin{lemme}\label{9.2}
Soit
\begin{equation*}
 \begin{CD}
  0@> >>\mathbb{Z} @> \imath >> \Gamma @>p>> G @> >> 1
 \end{CD} 
\end{equation*}
une extension centrale d'un groupe de Lie $G$ semi-simple {\`a} centre fini et
connexe. Notons $T=\imath(1)$. Il existe au plus une application $
\Phi : \Gamma \longrightarrow \mathbb{R}$ telle que
\begin{itemize}
\item[$(i)$] $\Phi$ est continue 
\item[$(ii)$] $\Phi(\gamma T) = \Phi(\gamma) +1$, $\forall \gamma \in \Gamma$
\item[$(iii)$] $ \Phi(\gamma_1 \gamma_2)-\Phi(\gamma_1)-\Phi(\gamma_2)$
est born{\'e}e sur $\Gamma\times \Gamma$
\item[$(iv)$] $\Phi(\gamma^n)=n\Phi(\gamma)$, $\forall \gamma \in
\Gamma$, $\forall n\in \mathbb{Z}$.
\end{itemize}
\end{lemme}
\proof Soit $\Phi_1$ et $\Phi_2$ deux telles fonctions, D'apr{\`e}s
$(ii)$, la fonction $\Phi_1-\Phi_2$ passe au quotient et d{\'e}finit une
fonction $F$ sur $G$. D'apr{\`e}s $(iii)$ elle est uniform{\'e}ment born{\'e}e sur les
commutateurs, et donc sur $N$ et sur $A$. D'apr{\`e}s $(i)$ elle est born{\'e}e
sur $K$. Utilisant la d{\'e}composition d'Iwasawa et {\`a} nouveau $(iii)$, on voit
que
$F$ est uniform{\'e}ment born{\'e}e sur $G$. D'apr{\`e}s $(iv)$, cela implique que
$F$ est identiquement nulle. D'o{\`u} $\Phi_1=\Phi_2$.
\hfill\hfill$\square$\\

Une telle fonction $\Phi$ est appel{\'e}e un {\it quasi-morphisme
homog{\`e}ne\/} \footnote{Notre d{\'e}finition diff{\`e}re l{\'e}g{\`e}rement de
celle de \cite{Barge-Ghys}.}.  Un  quasi-morphisme d{\'e}finit, par passage au quotient, une
application
$\Psi : G\longrightarrow {\mathbb{R}/\mathbb{Z}}$.\\

Revenant aux notations du paragraphe 5, on va maintenant construire un
quasi-morphisme pour le groupe $\Gamma$.

Soit
$\widetilde{o}$ un point base dans $\widetilde{S}$, et posons pour 
$\gamma\in \Gamma$
$$ c(\gamma) = m(\gamma\cdot\widetilde{o}, \widetilde{o})\ .$$
La primitive $m(\widetilde{\sigma}, \widetilde{\tau})$ est
invariante par $\Gamma$, et donc pour tous $\gamma_1,\gamma_2\in \Gamma$
$$\begin{array}{ll}
c(\gamma_1\gamma_2)- c(\gamma_1)-c(\gamma_2)&=
m(\gamma_1\gamma_2\cdot\widetilde{o}, \widetilde{o})-
m(\gamma_1\cdot\widetilde{o},\widetilde{o})-
m(\gamma_2\cdot\widetilde{o},\widetilde{o})\\
&=m(\gamma_1\gamma_2\cdot\widetilde{o}, \widetilde{o}) + 
m(\widetilde{o},\gamma_1\cdot\widetilde{o})+ 
m(\gamma_1\cdot\widetilde{o}, \gamma_1\gamma_2\cdot\widetilde{o})
\end{array}
$$
qui d'apr{\`e}s (\ref{Formule_Leray_Gen}) est l'indice de Maslov des projections sur $S$ des points
$\gamma_1\gamma_2\cdot\widetilde{o}$, $\widetilde{o}$ et
$\gamma_1\cdot\widetilde{o}$ de sorte que
\begin{equation*}
|c(\gamma_1\gamma_2)-
c(\gamma_1)-c(\gamma_2)|\leq r \ .
\end{equation*}
Pour tout {\'e}l{\'e}ment $\gamma \in \Gamma$, la suite $ (c_k):=(c(\gamma^k))_{k\in \mathbb{N}}$ satisfait donc la propri{\'e}t{\'e}
\begin{equation*}
\forall k,\, \ell\in \mathbb{N},\quad |c_{k+\ell}-c_k-c_\ell|\leq r\ .
\end{equation*}
D'apr{\`e}s un r{\'e}sultat classique, on d{\'e}duit que
\begin{equation*}
\tau(\gamma)=\lim_{k\rightarrow +\infty}{\frac{1}{k}}\,c(\gamma^k)
\end{equation*}
existe pour tout {\'e}l{\'e}ment $\gamma\in \Gamma$.

\begin{theoreme} 
La fonction  $$-\frac{1}{2}\tau$$ est un
quasi-morphisme sur
$\Gamma$. Elle est ind{\'e}pendante du choix de $\widetilde{o}$.
\end{theoreme}
\proof
La continuit{\'e} de $\tau$ se d{\'e}montre comme dans \cite[Proposition 2.11]{Barge-Ghys}.  Les
autres propri{\'e}t{\'e}s sont des cons{\'e}quences faciles des propri{\'e}t{\'e}s de la
primitive $m$ et notamment (\ref{action_T_m}). L'ind{\'e}pendance du point base
r{\'e}sulte du lemme \ref{9.2}, compte-tenu de ce que le groupe $G$ est
semi-simple, {\`a} centre fini et connexe.

\hfill\hfill$\square$\\

Par passage au quotient le quasi-morphisme $-\frac{1}{2}\tau$ d{\'e}fini une
fonction de $G$ dans $\mathbb{R}/\mathbb{Z}$ par 
$$\rho(g) =
-\frac{1}{2}\tau(\gamma)\, \mod\mathbb{Z}$$ o{\`u} $\gamma$ est un
rel{\`e}vement quelconque de $g$. Pour $g\in G$, la quantit{\'e} $\rho(g)$ est appel{\'e}e {\it nombre
  de rotation g{\'e}n{\'e}ralis{\'e}\/}  de $g$. La fonction
$\rho$ est invariante par conjugaison (toujours d'apr{\`e}s l'unicit{\'e}).

\begin{proposition}
La fonction $\rho$ poss{\`e}de les propri{\'e}t{\'e}s
suivantes :
\begin{itemize}
\item[$(i)$] si $g\in G$ fixe un point de $S$, alors $\rho(g)=0\ \mod \mathbb{Z}$.
\item[$(ii)$] si $u\in U$, alors $ e^{2i\pi\rho(u)}=\chi(u)$.
\end{itemize}
\end{proposition}
\proof
Pour $(i)$, on choisit comme origine $\widetilde{o}$ dans $\widetilde{S}$ un
point qui se projette sur un point de $S$ fix{\'e} par $g$, et on rel{\`e}ve $ g$ en un {\'e}l{\'e}ment
$\gamma$ tel que $\gamma\cdot\widetilde{o}=\widetilde{o}$, et donc aussi $\gamma^k\cdot\widetilde{o} = \widetilde{o}$ pour tout $k\in \mathbb{N}$. Il en r{\'e}sulte que
$\tau(\gamma)=0$.\\
Pour $(ii)$, choisissons un r{\'e}el $\varphi$ tel que $\chi(u) = e^{ir\varphi}$.
Un rel{\`e}vement de l'{\'e}l{\'e}ment $u$ est l'{\'e}l{\'e}ment
$\gamma$ dont l'action sur $\widetilde{S}$ est d{\'e}finie par
$$\gamma\cdot\widetilde{\sigma}=\gamma\cdot(\sigma,\theta)=(u(\sigma), \theta +\varphi),$$
puisqu'en effet $ \det(u(\sigma)) = \chi(u)\det(\sigma) =
e^{ir\varphi}e^{ir\theta} = e^{ir(\theta+\varphi)}.$
Donc$$\frac{1}{k} c(\gamma^k) = \frac{1}{k}\
m(\gamma^k\cdot\widetilde{\sigma},\widetilde{\sigma})=\frac{1}{k}\frac{1}{\pi}(\Psi(u^k(\sigma),
\sigma)-kr\varphi)
\rightarrow -\frac{{\varphi}r}{\pi}$$
quand $k\rightarrow +\infty$, puisque $\Psi$ est une fonction born{\'e}e.
Par suite,
$$ \rho(u)=\frac{\varphi r}{2\pi}$$
et donc $ e^{2i\pi\rho(u) } = e^{ir\varphi}=\chi(u)$.
\hfill\hfill$\square$

\end{document}